\PassOptionsToPackage{table}{xcolor} 
\documentclass[10pt]{amsart}
\usepackage[utf8x]{inputenc}
\usepackage[final]{changes}
\definechangesauthor[name={AR}, color=red]{ar}

\definecolor{lightgray}{gray}{0.9}
\definecolor{lightgreen}{rgb}{0.07, 0.04, 0.56} 
\definecolor{bluemunsell}{rgb}{0.0, 0.5, 0.69}
\definecolor{carolinablue}{rgb}{0.6, 0.73, 0.89}
\definecolor{turquoisegreen}{rgb}{0.63, 0.84, 0.71}
\definecolor{silver}{rgb}{0.75, 0.75, 0.75}
\definecolor{slategray}{rgb}{0.44, 0.5, 0.56}
\definecolor{blond}{rgb}{0.98, 0.94, 0.75}
\definecolor{lightgray}{gray}{0.9}

\usepackage{tikz}
\usepackage{tkz-tab}
\usepackage{tkz-graph}
\usepackage{pgfkeys}
\usetikzlibrary{shapes,arrows,fit,calc,positioning, decorations.pathmorphing, fit, petri}
\usetikzlibrary{arrows.meta, quotes}
\tikzstyle{vertex}=[circle, draw, inner sep=0.05pt, minimum size=0pt]

\tikzset{
  invisible/.style={opacity=0},
  visible on/.style={alt=#1{}{invisible}},
  alt/.code args={<#1>#2#3}{%
    \alt<#1>{\pgfkeysalso{#2}}{\pgfkeysalso{#3}} 
  },
  nodefill/.style={alt=<#1>{fill=blue!30,rounded corners}{},anchor=base}
}

\tikzset{anchor/.append code=\let\tikz@auto@anchor\relax}
\makeatother
\tikzset{
    clock hours/.style={
        label={ [ vertex, inner sep=5pt, double, nodefill=1]-(120*#1-150):$z_{{#1}}$}}
}
\tikzset{
   FiveEcksPlusOne/.style={
        label={ [ vertex, inner sep=5pt, double, nodefill=1]-(120*#1-150):$#1$}}
}

\tikzset{
   PolynomialCycle/.style={
        label={ [ vertex, inner sep=5pt, double, nodefill=1]-(90*#1-270):$z_{#1}$}}
}

\tikzset{
   PolynomialCycle3D/.style={
        label={ [ vertex, inner sep=2pt, double, nodefill=1]-(90*#1-270):\tiny{$z_{#1}$}}}
}

\usepackage{amsmath}
\usepackage{amsfonts}
\usepackage{amsthm}
\usepackage{amssymb}

\usepackage{setspace}
\usepackage[lined, boxed,linesnumbered, ruled]{algorithm2e}
\providecommand{\keywords}[1]{\textbf{\textit{Index terms---}} #1}
\usepackage{enumerate}
\usepackage{graphicx}
\usepackage{setspace}

\makeatletter
\newcommand\footnoteref[1]{\protected@xdef\@thefnmark{\ref{#1}}\@footnotemark}
\makeatother
\usepackage{caption, subcaption}
\usepackage{pdfpages}

\newtheorem{theorem}{Theorem}[section]
\newtheorem{hypotheses}[theorem]{Assumptions}
\newtheorem{lemma}[theorem]{Lemma}

\newtheorem{identities}[theorem]{Identities}
\newtheorem{corollary}[theorem]{Corollary}
\newtheorem{definition}[theorem]{Definition}

\usepackage{fancyhdr}

\title[{A Dual-Radix Modular Division Algorithm}]{A Dual-Radix Modular Division Algorithm for Computing Periodic Orbits within Syracuse Dynamical Systems}
\author{Andrey Rukhin}

\newcommand{\addend}[1]{\ensuremath{\addendSymbol_{#1}}}

\newcommand{\addendSymbol}{\ensuremath{r}}
\newcommand{\admissibleFactor}[1]{\ensuremath{\admissibleFactorSymbol_{#1}}}
\newcommand{\admissibleFactorTuple}[1]{\ensuremath{\iterateRemainder{#1},\translationValue{#1},\height{#1}}}
\newcommand{\admissibleFactorSymbol}{\ensuremath{t}}
\newcommand{\admissiblePair}[1]{\ensuremath{\left(\exponent{#1},\addend{#1}\right)}}
\newcommand{\admissibleResidueSet}[1]{\ensuremath{\mathcal{A}_{#1}}}

\newcommand{\bfrac}[2]{\left[\ensuremath{\frac{#1}{#2}}\right]}

\newcommand{\brackets}[1]{\ensuremath{\left[#1\right]}}
\newcommand{\bracketsInverse}[2][1]{\ensuremath{\brackets{#2}^{-#1}}}

\newcommand{\crossSection}[1]{\ensuremath{\graded{c}}_{#1}}

\newcommand{\denominator}[1]{\ensuremath{\denominatorSymbol_{#1}}}
\newcommand{\denominatorSymbol}{\ensuremath{D}}

\newcommand{\digitDifference}[1]{\ensuremath{\digitDifferenceSymbol_{#1}}}
\newcommand{\digitDifferenceSymbol}{\ensuremath{c}}

\newcommand{\dualRadixMappings}{\ensuremath{\mathcal{D}}}

\newcommand{\equivMod}[1]{\ensuremath{\underset{#1}{\equiv}}}

\newcommand{\exponent}[1]{\ensuremath{\exponentSymbol_{#1}}}

\newcommand{\exponentSymbol}{\ladicExponentSymbol}

\newcommand{\exponentSum}[1]{\less{#1}}

\newcommand{\factor}[1]{\ensuremath{{\bf #1}}}

\newcommand{\forwardExponent}[1]{\ensuremath{\forwardExponentSymbol_{#1}}}
\newcommand{\forwardExponentFactor}{\ensuremath{\factor{\forwardExponentSymbol}}}

\newcommand{\forwardExponentSum}[1]{\meps{#1}}
\newcommand{\forwardExponentSymbol}{\madicExponentSymbol}

\newcommand{\forwardMapSymbol}{\ensuremath{\mathcal{L}}}
\newcommand{\forwardMapFullSymbol}[1]{\ensuremath{\forwardMapSymbol_{#1}}}
\newcommand{\forwardMap}[2]{\ensuremath{\forwardMapSymbol_{#2}\left(#1\right)}}
\newcommand{\forwardMapFull}[3]{\ensuremath{\forwardMapFullSymbol{#1; #3}\left(#2\right)}}

\newcommand{\gradationSymbol}{\ensuremath{g}}
\newcommand{\graded}[1]{\ensuremath{{\bf #1}}}

\newcommand{\height}[1]{\ensuremath{\heightSymbol_{#1}}}

\newcommand{\heightSymbol}{\ensuremath{h}}

\newcommand{\integers}{\ensuremath{\mathbb{Z}}}
\newcommand{\integersMod}[1]{\ensuremath{\mathbb{Z}}/{#1\integers}}

\newcommand{\inverseMapSymbol}{\ensuremath{\mathcal{D}}}
\newcommand{\inverseMapFullSymbol}[1]{\inverseMapSymbol_{#1}}
\newcommand{\inverseMapFull}[4]{\ensuremath{\inverseMapFullSymbol{#1}\left(#2;#3; #4\right)}}

\newcommand{\iterate}[1]{\ensuremath{\iterateSymbol_{#1}}}
\newcommand{\iterateRemainder}[1]{\ensuremath{\iterateRemainderSymbol_{#1}}}
\newcommand{\iterateRemainderSymbol}{\ensuremath{s}}
\newcommand{\iterateSymbol}{\ensuremath{n}}

\newcommand{\ladicDigit}[1]{\ensuremath{\ladicDigitSymbol_{#1}}}
\newcommand{\ladicDigitSymbol}{\ensuremath{b}}

\newcommand{\ladicExponentSymbol}{\ensuremath{e}}
\newcommand{\lExp}[1]{\ladicExponentSymbol_{#1}}
\newcommand{\leps}[1]{\ladicExponentPrefixSum{#1}}
\newcommand{\less}[1]{\ladicExponentSuffixSum{#1}}
\newcommand{\ladicExponentPrefixSum}[1]{\ensuremath{\ladicExponentSumSymbol_{#1}}}
\newcommand{\ladicExponentSuffixSum}[1]{\ensuremath{\overline{\ladicExponentSumSymbol}_{#1}}}

\newcommand{\ladicExponentSumSymbol}{\ensuremath{E}}

\newcommand{\ladicPrefixAddend}[1]{\ensuremath{\ladicPrefixAddendSymbol_{#1}}}
\newcommand{\ladicPrefixAddendSymbol}{\ensuremath{Q}}
\newcommand{\ladicQuotient}[1]{\ensuremath{\ladicQuotientSymbol_{#1}}}
\newcommand{\ladicQuotientSymbol}{\ensuremath{j}}

\newcommand{\ladicResidue}[1]{\ensuremath{\ladicResidueSymbol_{#1}}}
\newcommand{\ladicResidueSymbol}{\ensuremath{\lambda}}
\newcommand{\ladicShift}[1]{\ensuremath{\ladicQuotient{#1}}}

\newcommand{\leftShift}{\langle}
\newcommand{\leftShiftOperator}[1]{\leftShift #1}

\newcommand{\List}[3]{\ensuremath{\left(#1_{#2},\ldots,#1_{#3}\right)}}

\newcommand{\madicDigit}[1]{\ensuremath{\madicDigitSymbol_{#1}}}
\newcommand{\madicDigitSymbol}{\ensuremath{d}}

\newcommand{\madicExponentSymbol}{\ensuremath{f}}
\newcommand{\mExp}[1]{\madicExponentSymbol_{#1}}
\newcommand{\meps}[1]{\madicExponentPrefixSum{#1}}
\newcommand{\mess}[1]{\madicExponentSuffixSum{#1}}

\newcommand{\madicExponentPrefixSum}[1]{\ensuremath{\madicExponentSumSymbol_{#1}}}
\newcommand{\madicExponentSuffixSum}[1]{\ensuremath{\overline{\madicExponentSumSymbol}_{#1}}}
\newcommand{\madicExponentSumSymbol}{\ensuremath{F}}

\newcommand{\madicPrefixAddend}[1]{\ensuremath{\madicPrefixAddendSymbol_{#1}}}
\newcommand{\madicPrefixAddendSymbol}{\ensuremath{P}}
\newcommand{\madicQuotient}[1]{\ensuremath{\madicQuotientSymbol_{#1}}}
\newcommand{\madicQuotientSymbol}{\ensuremath{k}}

\newcommand{\madicResidue}[1]{\ensuremath{\madicResidueSymbol_{#1}}}
\newcommand{\madicResidueSymbol}{\ensuremath{\mu}}
\newcommand{\madicShift}[1]{\ensuremath{\madicQuotient{#1}}}

\newcommand{\mixedRadixAdics}[2]{Z_{\factor{#1}, \factor{#2}}}
\newcommand{\mlSystem}{\ensuremath{\left(m,l, \translationMatrix, \forwardExponentFactor \right)}}

\newcommand{\naturals}{\ensuremath{\mathbb{N}}}
\newcommand{\nequivMod}[1]{\ensuremath{\underset{#1}{\not\equiv}}}
\newcommand{\nOpenSet}[1]{\ensuremath{\left[#1\right)}}

\newcommand{\numerator}[1]{\ensuremath{\numeratorSymbol_{#1}}}
\newcommand{\numeratorBase}[1]{\ensuremath{\numeratorBaseSymbol_{#1}}}
\newcommand{\numeratorBaseSymbol}{\ensuremath{\gamma}}
\newcommand{\numeratorSymbol}{\ensuremath{N}}
\newcommand{\nZeroOpenSet}[1]{\ensuremath{\left[#1\right)_0}}

\newcommand{\nZeroSet}[1]{\ensuremath{\left[#1\right]_0}}

\newcommand{\ord}[2]{\ensuremath{\nu_{#1}\left(#2\right)}}

\newcommand{\padics}[1]{\ensuremath{Z_{#1}}}

\newcommand{\parentheses}[1]{\ensuremath{\left(#1\right)}}

\newcommand{\pdiv}{\mid\!\mid}
\newcommand{\pfrac}[2]{\ensuremath{\left(\frac{#1}{#2}\right)}}

\newcommand{\prefixAddend}[1]{\ensuremath{\prefixAddendSymbol_{#1}}}
\newcommand{\prefixAddendSymbol}{\ensuremath{P}}

\newcommand{\rationals}{\ensuremath{\mathbb{Q}}}
\newcommand{\rationalUnits}[1]{\ensuremath{\mathbb{Z}_{\left\langle#1\right\rangle}}}

\newcommand{\resRep}[2]{\ensuremath{\brackets{#1}_{#2}}}

\newcommand{\rightShift}{\rangle}
\newcommand{\rightShiftOperator}[1]{\rightShift #1}

\newcommand{\translationMatrix}{\ensuremath{\factor{A}}}
\newcommand{\translationMatrixEntry}[2]{\ensuremath{\translationMatrixEntrySymbol_{#1,#2}}}
\newcommand{\translationMatrixEntrySymbol}{\ensuremath{a}}
\newcommand{\translationValue}[1]{\ensuremath{\translationMatrixEntrySymbol_{#1}}}
\newcommand{\translationValueSymbol}{\ensuremath{\translationMatrixEntrySymbol}}

\renewcommand{\emph}[1]{{\sl #1}}

\begin{document}

\begin{abstract}

This article analyzes the periodic orbits of  Syracuse dynamical systems in a novel algebraic setting: the commutative ring of {\sl graded $n$-adic integers}.   Within this context, this article introduces a {\sl dual-radix} modular division algorithm for computing the {\sl graded canonical expansions} and {\sl graded quotients} for a certain class of rational expressions that arise from periodic orbits within these dynamical systems.   This division algorithm yields two novel methods for testing the integrality of the B\"{o}hm-Sontacchi numbers \cite{BohmSontacchi78}.

 \end{abstract}

\keywords{
$3x+1$ Problem, Collatz Conjecture, Kakutani's Problem, Syracuse Algorithm, Parallel $p$-adic Division
}

\maketitle

\section{Introduction}

	\subsection{Background and Motivation}

	 Let  $m,l$ be coprime integers exceeding $1$, and let $f \in \naturals$. In this article, we will study a family of dynamical systems derived from the iterations of mapping  on $\mathbb{Q}$ of the form
\[
x \to \frac{m^{\forwardExponentSymbol}x + \translationMatrixEntrySymbol}{l^{\exponentSymbol}}
\] where $a \ (a\in \integers)$ is such that  $e = \ord{l}{m^fx+a} > 0$.   Modeled after the \emph{(accelerated) $3x+1$ dynamical system}, this class of dynamical systems  encompasses its progenitor;  other generalizations can be found in \cite{KontorovichSinai}, \cite{Conway},  and \cite{MatthewsWatts}.

 One open research area considers the analyses of periodic orbits within such dynamical systems.  
 It has been shown (within the context of the (accelerated) $3x+1$ Problem) that the iterate values  within a periodic orbit of length $\tau$ assume a rational expression of the form
 \begin{equation}\label{equation::RationalExpression}
	\frac{\numerator{}}{\denominator{}} = \frac{\sum_{0 \leq w < \tau}m^{\forwardExponentSum{w}}l^{\exponentSum{\tau-1-w}}\translationMatrixEntrySymbol_w }{l^{\exponentSum{\tau}} - m^{\forwardExponentSum{\tau}}}\hspace{2pt}
	\end{equation} 
	where the sequences $\left(\meps{0},\ldots,\meps{\tau}\right)$ and $\left(\less{0},\ldots,\less{\tau}\right)$ are strictly increasing over $\naturals_0$  with $\meps{0} = \less{0} = 0$, and the term  $\translationValue{w} \in \integers$ for $w\in \{0,\ldots,\tau-1\}$.
	As an example, when the base $m=3$, the base $l=2$,   the exponent $\forwardExponentSum{w} = w$, and  the term $\translationValue{w} = 1$ for $w\in \{0,\ldots,\tau-1\}$, these rational expressions yield the iterate values in an (accelerated) orbit of length $\tau$ within the $3x+1$ dynamical system (of unaccelerated length $ \less{\tau}$); in this context, such  expressions  are known as the {\sl B\"{o}hm-Sontacchi numbers} \cite{BohmSontacchi78}.  \added{The mixed-base representations of the B\"{o}hm-Sontacchi numbers are considered in   \cite{Anisimov2009},  they have been studied in further detail (within the context of rationals with odd denominators) in \cite{Lagarias1990}, and the {\sl 3-smooth representations} of the numerators have been studied in \cite{Avidon} and \cite{Selfridge}.}  A notorious open question pertains to the existence of positive integers (exceeding $1$) that admit such rational expressions (see  \cite{Lagarias10} and  \cite{wirsching2006dynamical} for  comprehensive surveys on this subject). 

 In this article, we will call such rational expressions the {\sl dual-radix rational expressions} of the iterate values. Both the $m$-adic and $l$-adic expansions of dual-radix rationals are well-defined, as the denominator $l^{\exponentSum{\tau}} - m^{\forwardExponentSum{\tau}}$  is coprime to both $m$ and $l$ (by virtue of $m$ and $l$ being coprime).  Herein, we will consider these expressions within the rings of {\sl graded $m$-adic} and {\sl graded $l$-adic integers} (introduced below).  This article presents a new modular division algorithm for computing the graded $m$-adic and graded $l$-adic expansions of such rational expressions; some features of this algorithm are:
\begin{itemize}
	\item (Theorem \ref{theorem::quotientDigitRecurrence}) for each $u\in \naturals$, the algorithm computes the $u$-th graded $m$-adic and  graded $l$-adic digits of the expansions of such rational expressions by way of the arithmetic difference of the $(u-1)$-th  graded $m$-adic and  graded $l$-adic digits (without carry-propagation, and without intermediate values that grow proportionally with $u$);
	\item (Theorem \ref{theorem::quotientDigitRecurrence}) for any desired mantissa length $\mu$, the algorithm compresses (by half) $2\tau$ \added{graded} Hensel codes of length $\mu$ (the codes of the graded $m$-adic and graded $l$-adic expansions for each of the $\tau$ orbit elements) using $\mu\tau$ values (each value in the compressed representation is an arithmetic difference of graded $m$-adic and graded $l$-adic digits);
		\item  (Theorem \ref{theorem::quotientDigitRecurrence}) the algorithm  computes the $m$-adic and $l$-adic canonical expansions of all $\tau$ iterate values in parallel;
	\item (Theorem \ref{theorem::mixedRadixQAdicQuotientShift}) these arithmetic differences of digits  in this compressed representation can be used to  express the graded $m$-adic and graded $l$-adic quotients (or {\sl shifts}\/) of the iterate values as a dual-radix rational.
\end{itemize}

For example, consider the following periodic orbit within the (accelerated) $3x-1$ dynamical system:
	\[
	17  \xrightarrow{\frac{3(17)-1}{2^1}} 25 \xrightarrow{\frac{3(25)-1}{2^1}} 37 \xrightarrow{\frac{3(37)-1}{2^1}} 55 \xrightarrow{\frac{3(55)-1}{2^2}} 41 \xrightarrow{\frac{3(41)-1}{2^1}} 61\xrightarrow{\frac{3(61)-1}{2^1}} 91 \xrightarrow{\frac{3(91)-1}{2^4}} 17 \to \cdots.
	\]

\begin{figure}[t]
		\vspace{70pt}
		\begin{subfigure}[b]{.45\linewidth}\label{Figure:3DExampleI}
				\hspace{-50pt}
				\scalebox{.70}{
					\includegraphics[scale=1,bb=20 20 120 120]{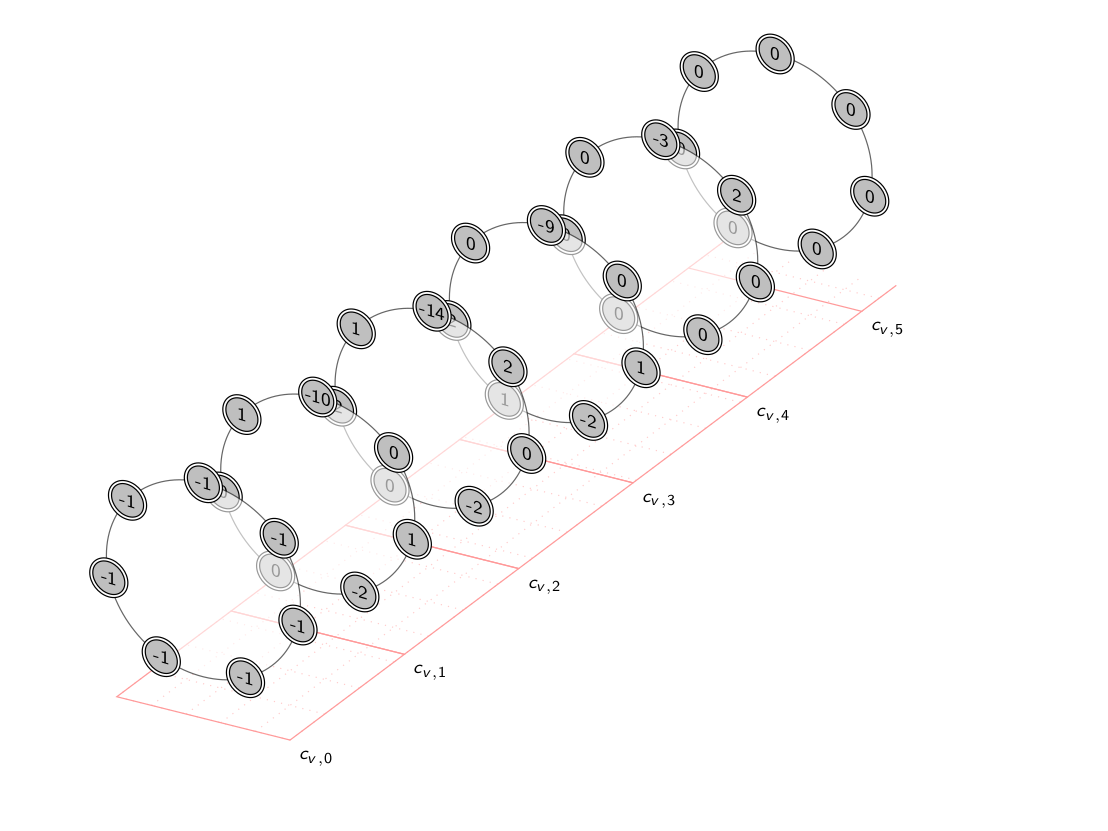}
				}
				\subcaption{A three-dimensional representation of the first six cross-sections of the quotient cylinder of the periodic orbit $(17, 25, 37, 55, 41, 61, 91)$ within the (accelerated) $3x-1$ dynamical system (Theorem \ref{theorem::quotientDigitRecurrence})}
				
		\end{subfigure} \hspace{20pt}
		\begin{subtable}[b]{.45\linewidth}
			\centering
\renewcommand{\arraystretch}{1.5}

	\begin{tabular}{|c||c|c|c|c|c|c|} \hline 
		$u\in \naturals_0 $ & $0$ & $1$ & $2$ & $3$ & $4$ & $5$ \\ \hline \hline
		$\digitDifference{0,u}$ & -1 & -10 & -14 & -9 & -3 & 0 \\ \hline
		$\digitDifference{1,u}$ & -1 & 0 & 2 & 0 & 2 & 0 \\ \hline
		$\digitDifference{2,u}$ & -1 & 1 & 0 & 1 & 0 & 0 \\ \hline
		$\digitDifference{3,u}$ & -1 & -2 & -2 & -2 & 0 & 0 \\ \hline
		$\digitDifference{4,u}$ & -1 & 0 & 0 & 1 & 0 & 0 \\ \hline
		$\digitDifference{5,u}$ & -1 & 0 & 2 & 2 & 0 & 0 \\ \hline
		$\digitDifference{6,u}$ & -1 & 1 & 1 & 0 & 0 & 0 \\ \hline
	\end{tabular}

			\subcaption{A two-dimensional representation of these cross-sections}
				\label{Figure:3DExampleI}
		\end{subtable}
	\caption{}\label{Figure:ExampleI}
	\end{figure}

	For this orbit, the dual-radix algorithm yields the arrangement of integers--the {\sl quotient cylinder}--illustrated in Figure \ref{Figure:ExampleI}. This arrangement compresses the $3$-adic and graded $2$-adic expansions of the seven periodic iterates (Theorem \ref{theorem::quotientDigitRecurrence}); furthermore, one can compute the $3$-adic and graded $2$-adic quotients of these iterates as dual-radix rationals with coefficient weights found along the cross-sections (or columns) of this arrangement (Theorem \ref{theorem::mixedRadixQAdicQuotientShift}).
	
Theorem  \ref{theorem::mixedRadixQAdicQuotientShift} also yields a  method (the {\sl suffix test}) for deciding the integrality of the rational expression in (\ref{equation::RationalExpression}); for the $3x+d$ Problem ($d = \pm 1$), this test can be performed at the $\tau$-th stage of the algorithm described herein by appealing to an upper-bound on the iterates established in \cite{BelagaMignotte}, along with upper bounds on linear forms of two logarithms   established in \cite{SimonsDeWeger} (based on the results in \cite{BakerWustholz} and \cite{Rhin}).

 Theorem \ref{theorem::NumeratorAsDifference} yields a method (the {\sl prefix test}) for deciding the integrality of the rational expression in (\ref{equation::RationalExpression}) by taking into account the prefixes of the $m$-adic and $l$-adic expansions of $\frac{N}{D}$.   Let $\madicResidue{\tau}\in \{0,\ldots,m^{\meps{\tau}}-1\}$ be such that $\madicResidue{\tau}\equiv ND^{-1} \bmod m^{\meps{\tau}}$, and let  $\ladicResidue{\tau}\in \{0,\ldots,l^{\less{\tau}}-1\}$ be such that $\ladicResidue{\tau}\equiv ND^{-1}\bmod {l^{\less{\tau}}}$. 
We will establish the identities
	\[
	\frac{\sum_{0\leq w < \tau}m^{\meps{w}}l^{\less{\tau-1-w}}\translationValue{w}}{l^{\less{\tau}}-m^{\meps{\tau}}} =  \madicResidue{\tau} + m^{\meps{\tau}}\pfrac{\madicResidue{\tau} - \ladicResidue{\tau}}{l^{\less{\tau}}-m^{\meps{\tau}}} = \ladicResidue{\tau} + l^{\less{\tau}}\pfrac{\madicResidue{\tau} - \ladicResidue{\tau}}{l^{\less{\tau}}-m^{\meps{\tau}}}
	\] in Theorem \ref{theorem::NumeratorAsDifference}.	
Consequently, a {\sl $3$-smooth representation}\\
of a natural number $n = \sum_{0\leq w < \tau}3^{\meps{w}}2^{\less{\tau-1-w}}$  (see \cite{Erdos})  can be expressed as the difference
\[
n = 2^{\less{\tau}}\resRep{nD^{-1}}{3^{\meps{\tau}}} - 3^{{\meps{\tau}}}\resRep{nD^{-1}}{2^{\less{\tau}}}
\] where $D =  2^{\less{\tau}} - 3^{\meps{\tau}}$,
  and the difference $D$ divides $n$ if and only if $D$ divides the difference of canonical representatives of residues $\resRep{nD^{-1}}{3^{{\meps{\tau}}}} - \resRep{nD^{-1}}{2^{\less{\tau}}}$.

 These methods will be applied in a future article to produce alternative proofs of Steiner's $1$-Cycle Theorem (\cite{Steiner}) within the context of the $3x+1$ Problem;   furthermore, these approaches show that $(1)$ and $(5,7)$ are the only $1$-cycles (over \naturals) within the (accelerated) $3x-1$ dynamical system for all but finitely many $\tau$.

  \subsection{Article Summary}
  
  We conclude Section 1 by detailing the notation used throughout this article, and we will introduce the commutative ring of {\sl graded $n$-adic integers}.  In  Section \ref{section::MLSystem}, we will detail the class of dynamical systems that serve as the framework of this article.
  In Section \ref{section::OrbitTerms}, we expand upon the result of  B\"{o}hm and Sontacchi by identifying the dual-radix rational expressions of periodic orbit elements within our dynamical systems;  furthermore, we will establish the relevant terminology for Section \ref{section::MainResults}, where  the primary results of this article are presented.  \deleted{In this final section, we will demonstrate a novel recurrence for computing the graded $m$-adic and graded $l$-adic digits of all $\tau$ orbit elements simultaneously, and we will also show how this recurrence yields a method for efficiently computing the graded quotients of the orbit elements. 	}

	We now proceed by  defining the notation used throughout this article.
\subsection{Notation Summary}
 If $a,b$ and $n$ are integers with $n\neq 0$, we  write $a\equivMod{n}b$ if
and only if $a$ is equivalent to $b$ modulo $n$.  

Context permitting\footnote{When it is well-defined, we will also write $\bracketsInverse{a} \bmod b$ (e.g., in Definition \ref{definition::admissibleResidues}) to denote the multiplicative inverse of $a$ in  $\integers/ b\integers$.}, we  denote the set $\{1,\ldots,i\}$ of the first $i$ positive integers as $[i]$.  We  write $[i)$ to denote the set $[i-1]$. 
We
 also write $[i]_0$ and $[i)_0$ to denote the sets $[i]\cup \{0\}$ and  $[i)\cup \{0\}$, respectively.

For  $b\in \naturals$  where $b\geq 2$ and $n\in \mathbb{Z}$,  we assume the standard definition of $b$-adic order: we will write

\begin{equation}
\ord{b}{n} = 
	\begin{cases}
		\max\left\lbrace u\in \naturals_0 : b^u \mid n \right\rbrace & n\neq 0 \\
		0 & n=0.
	\end{cases}\notag
\end{equation}
  For a rational $\frac{n}{d} \in \rationals$, we assume the definition $\ord{b}{\frac{n}{d}} = \ord{b}{n} - \ord{b}{d}$.  We may also write $b^u \pdiv n$ if and only if $\ord{b}{n} = u$.

		We adopt the convention of denoting sequences of elements over a set with bold notation.  Let $\factor{x} = \List{x}{0}{\tau-1}$ be a sequence of length $\tau$ over a set $X$, and let $v\in \integers$.  We will denote the \emph{left cyclic shift (of index $v$) modulo $\tau$} of  \graded{x} as $\graded{x} {\leftShift} {v}$, where
		$
		\replaced{\graded{x} \leftShift v = \List{x}{ v \bmod \tau}{\tau-1+v\bmod \tau} }{\graded{x}{\leftShift} v = \List{x}{ v \bmod z}{v + \tau-1\bmod z}}.
		$ We will also denote the  \emph{right cyclic shift (of index $v$) \added{modulo $\tau$}} of  \graded{x} as $\graded{x} {\rightShift} {v}$, where
		$\replaced{\graded{x}{\rightShift} v = \List{x}{ -v \bmod \tau}{\tau-1-v\bmod \tau} }{
		\graded{x}  {\rightShift} v= \List{x}{ -v \bmod z}{-v + z-1 \bmod z}}.
		$
		The shift notation may be omitted when the shift index equals zero.  
		
		
		We will let $\graded{x}^R$ denote the {\sl reversal} of the sequence $\graded{x}$.
		From the definitions, it follows that
		$
		\graded{x}^R \leftShiftOperator{v} = \left[\graded{x} \rightShift v\right]^R,
		$ and
		$
		\graded{x}^R \rightShiftOperator{v} = \left[\graded{x} \leftShift v\right]^R.
		$

		For any sequence $\graded{y} = \parentheses{y_w}_{w=0}^{\tau-1}$ of length $\tau$, we assume the standard definition of the dot product $\graded{x}\cdot \graded{y} = \sum_{0\leq w < \tau}x_wy_w$.	
		

		We will  often evaluate a function at shifted permutations of its arguments.   For $u\in \naturals$, let $\graded{x}_u = \left(x_{w\bmod \tau}\right)_{w=0}^{u-1}$,  and let $\phi_u$ be a function defined on $X^u$. We will adopt the above notational convention and write $\phi_{v,u}\left(\graded{x}_u\right)  := \phi_u\left(\graded{x}_u{\leftShift v}\right).$  


	We will now detail the \emph{graded $n$-adic} number system that  will be studied in this article.
\subsection{The Commutative Ring of Graded $n$-adic Integers} 

Let $l,m\in \integers$ where $l,m\geq 2$. We will define the multiplicative set 
	\[
	S_{l,m} = \left\lbrace s \ \middle| \    s \in \mathbb{Z},\   \gcd(s,l) = \gcd(s,m)=1 \right\rbrace.
	\]
	We will  write $\rationalUnits{l,m}$ to denote the localization of $\mathbb{Z}$ by $S_{l,m}$, and when $l=m$, we will   write $\rationalUnits{l}$.

The following construction of the commutative ring of graded $n$-adic integers is analogous to the algebraic construction of the standard $n$-adic integers (e.g., see Sutherland \cite{Sutherland}).  
Let $\tau$ be a positive integer.  Let \factor{g}  be  a sequence of positive integers of length $\tau$ where $\factor{g} = \List{g}{0}{\tau-1}$, and let $g_u = g_{u \bmod \tau}$ for $u\in \naturals_0$. For each $y\in \naturals_0$, define $G_y$ to be the prefix sum $\sum_{0\leq w < y}g_{w}$.

Let $R_u$ denote the ring $\integers/ n^{G_{u}}\integers$ for each $u\in \naturals$.  For each element $x$ of $R_u$, we adopt the graded representation 
$
x = \sum_{0\leq w < u} \madicDigit{w} n^{G_w}
$ where $\madicDigit{w} \in\nZeroOpenSet{n^{g_w}}$.

 Define the morphism $f_u: R_{u+1}\to R_u$ to be reduction modulo $n^{G_u}$.  We define the commutative ring of ({\sl graded}\/)\emph{ \factor{n}-adic integers} $Z_{\factor{n},\factor{g}}$ to be the inverse limit\footnote{The sequence $\left(G_u\right)_{u\geq 1}$ is cofinal in \naturals; thus, the inverse limit  $Z_{\factor{n},\factor{g}}$ is isomorphic to the ring of $n$-adic integers (see Ribes and Zalesskii \cite[p.8]{ribes2000profinite}). However, we will write \factor{n},\factor{g} instead of $n$ to emphasize the graded representation of the ring elements.}
\[
Z_{\factor{n},\factor{g}} = \varprojlim R_u
\] of the inverse system $(R_u, f_u)_{u\geq 1}$.
We define the sequence $\factor{g}^{\infty} = \left(g_u\right)_{u\in \naturals_0}$ to be the \emph{gradation} of $Z_{\factor{n},\factor{g}}$, where the sequence $\factor{g}$ is the  \emph{gradation sequence (or grading)} of $Z_{\factor{n},\factor{g}}$.  

 Similar to the standard $n$-adics, an  element  $\frac{a}{b} \in \rationalUnits{n}$  is mapped into $\padics{\factor{n},\factor{g}}$ by a ring homomorphism to the inverse limit element $\left(\prefixAddend{u}\right)_{u\geq 1}$, where $\prefixAddend{u}\in \nZeroOpenSet{n^{G_u}}$ and $\prefixAddend{u} \equiv ab^{-1} \bmod n^{G_u}.$  

\subsubsection{Graded Canonical Expansions}
One can write 
$
\prefixAddend{u} = \sum_{0\leq w < u} \madicDigit{w} n^{G_w},
$ where $\madicDigit{w} \in \nZeroOpenSet{n^{g_w}}$ for each $u\in \naturals$. 
Thus, we can uniquely express the element $\frac{a}{b}$ within $Z_{\factor{n},\factor{g}}$ with the series expansion
$
\sum_{w\geq 0} \madicDigit{w} n^{G_w}
$ where $\madicDigit{w}$ is the ({\sl graded}\/) \emph{$\factor{n}$-adic digit of $\frac{a}{b}$ of index $w$} (\emph{in, or over, the ring $Z_{\factor{n},\factor{g}}$}\/). 

As with the standard representation of the $n$-adics, we will define the sequence of digits $\left(\madicDigit{w}\right)_{w\geq 0}$ to be the \emph{canonical expansion of $\frac{a}{b}$} (\emph{in, or over, the ring $Z_{\factor{n},\factor{g}}$}\/).

\subsubsection{Graded Quotients}

Let $u \in \naturals$.  If we write
$
\frac{a}{b} = n^{G_u}\madicShift{u} + P_u
$ where $P_u \in \nZeroOpenSet{n^{G_u}}$ and $\madicShift{u}\in\rationalUnits{n}$,  then we define the quantity $\madicShift{u}$ to be the ({\sl graded}\/) {\sl \factor{n}-adic quotient (or shift)}\/  \emph{of index $u$} \emph{of $\frac{a}{b}$} (\emph{in, or over, the ring  $\mixedRadixAdics{n}{g}$}\/). We will  define the quantity $\prefixAddend{0}$ to be $0$, and we will define the quotient of index $0$ to be $\frac{a}{b}$.

\section{The \mlSystem\ Linear Dynamical System}\label{section::MLSystem}

To initiate our study of rational expressions in (\ref{equation::RationalExpression}), we  choose integers $m$ and $l$ satsifying Hypotheses \ref{hypotheses::mlAf} below, and
we will define a family of dynamical systems, called $(m,l)$-systems.  These  discrete dynamical systems are designed for analyzing periodic orbits  of a prescribed length within a dynamical system that scales  an argument $x$ (where $\ord{l}{x} = 0$)  by one or more powers of $m$, additively translates this scaled value to a value $x'$ (where $\ord{l}{x'}>0$), and then  contracts the translated value  by one or more powers of $l$ to a value $x''$ (where $\ord{l}{x''}=0$).     We do not require that our choices of scaling powers and additive translations are the same at each iteration.

\begin{definition}[$(m,l)$-Systems] 
Let $m$ and $l$ be coprime integers where  $m\geq 2$ and $l\geq 2$.  Fix $\tau\in \naturals$, and let $\factor{f}$ be a sequence of $\tau$ positive integers where $\factor{f} = \left(\forwardExponent{v}\right)_{v=0}^{\tau-1}$ (the \emph{$m$-adic grading} of the system).   Let $\translationMatrix  \in M(\integers)_{\tau,l}$ be a matrix with $\tau$ rows and $l$ columns indexed by the sets $\nZeroOpenSet{\tau}$ and $\nZeroOpenSet{l}$, respectively, in increasing order; furthermore, the entries of \translationMatrix\ satisfy the following conditions:
 for all $v\in \nZeroOpenSet{\tau}$,  if $i=0$, then the entry $\translationMatrixEntry{v}{i} = 0$, and  if $i\neq 0$, then:
\begin{enumerate}[i.] 
	\item  $\translationMatrixEntry{v}{i} \equivMod{l} -m^{\forwardExponent{v}}i$, 
	\item $\translationMatrixEntry{v}{i} \nequivMod{m}0$ and  $\translationMatrixEntry{v}{i} \nequivMod{l}0$. 
\end{enumerate}
Let $x \in \rationalUnits{l}$  where $x = \frac{y}{z}$, and $yz^{-1}\equivMod{l}i$ for some $i\in \nZeroOpenSet{l}$. For $v\in\nZeroOpenSet{\tau}$, define 
	$$\exponent{v} = \begin{cases}
						\ord{l}{m^{\forwardExponent{v}}x + \translationMatrixEntry{v}{i}} & i \neq 0,\\ 
						\ord{l}{x} & i = 0;
						 \end{cases}$$ 
We will define the \emph{\graded{l}-adic map} $\forwardMapFullSymbol{m,l,\translationMatrix; v}:\rationalUnits{l} \to  \rationalUnits{l}$ \emph{ of index $v$}  on $x$ to be
\[
\forwardMapFull{m,l,\translationMatrix}{x}{v} = 
	\begin{cases}
	\parentheses{m^{\forwardExponent{v}}x + \translationMatrixEntry{v}{i} }/ { l^{\exponent{v} } } & i\neq 0\\
	x/{ l^{\exponent{v} } } & i=0.
	\end{cases}
\]
We denote this ordered collection of $\tau$ \graded{l}-adic maps with the tuple $\mlSystem$, and we will define it to be a \emph{$(m,l)$-system} (\emph{of order $\tau$}).
\end{definition}

Context permitting, we will forgo writing the subscript(s) and  write $\forwardMapSymbol$.

Note that the exponent $\exponent{v}$ is positive due to our restrictions on the entries of the translation matrix $\translationMatrix$.

We will now establish the setting used throughout this article.

\begin{hypotheses}\label{hypotheses::mlAf}
Let $m$ and $l$ be coprime, positive integers where  $m\geq 2$ and $l\geq 2$.  Let $\tau$ be a positive integer.  Let $\factor{f}$ be a sequence of $\tau$ positive integers where $\factor{f} = \left(\forwardExponent{v}\right)_{v=0}^{\tau-1}$, and let  $\forwardExponent{u} =\forwardExponent{u\bmod \tau}$ for $u\in \integers$.  Finally, let $\translationMatrix \in M(\integers)_{\tau,l}$ be such that $\mlSystem$ is a  $(m,l)$-system of order $\tau$.
\end{hypotheses}

Our analysis of periodic orbits in a given \replaced{$(m,l)$-system}{dynamical system} necessitates determining \replaced{the inverse of an iterate value under an \graded{l}-adic mapping when it exists.}{the {\sl dual-radix mappings} of the system (defined below).} To this end,  the next definition identifies the possible residue classes mod  $m^{\forwardExponent{v}}$ that admit an invertible image \deleted{$\frac{m^{\forwardExponent{v}}x + \translationMatrixEntry{i}{v} }{ l^{\exponent{v} } }$}of an \graded{l}-adic mapping.

\begin{definition}[Admissible Residues]\label{definition::admissibleResidues}
 Assume \ref{hypotheses::mlAf}.

Let $i\in \nOpenSet{l}$, and let $v\in \nZeroOpenSet{\tau}$.  Let $\iterateRemainder{v} \in\nOpenSet{m^{\forwardExponent{v} }}$ where $\gcd(\iterateRemainder{v}, m)=1$. We will say that the ordered pair $\left(\iterateRemainder{v},\translationMatrixEntry{v}{i}\right)$ (where  $\translationMatrixEntry{v}{i}$ is an entry in $\translationMatrix$) is \emph{admissible} if and only if there exists an ordered pair  of positive integers $(\exponentSymbol,\addendSymbol)$  where  $\addendSymbol\equivMod{l} i$, and 
the pairs satisfy the equality 
		\begin{equation}
			\iterateRemainder{v} = \frac{m^{\forwardExponent{v}}\addendSymbol + \translationMatrixEntry{v}{i} }{l^{\exponentSymbol}};\label{equation::admissibilityEquation}
		\end{equation}  furthermore,  we   require that
	 $| \translationMatrixEntry{v}{i} |\in \nZeroOpenSet{\max\left(m^{\forwardExponent{v}}, l^{\exponentSymbol}\right)}$.

	 We will denote the set of admissible pairs  (of index $v$) as $\admissibleResidueSet{v}$. For a given  admissible pair $\left(\iterateRemainder{v}, \translationMatrixEntry{v}{i}\right)$, we will define such a $\left(\exponentSymbol, \addendSymbol\right)$ pair to be a \emph{$\left(\iterateRemainder{v},\translationMatrixEntry{v}{i}\right)$-witness}.

\end{definition}

The existence of admissible pairs is assured as follows: for each $e \in\naturals$, the $m^f$-residue $s$ in an admissible pair (\iterateRemainderSymbol, \translationMatrixEntrySymbol) is such that 
$
\iterateRemainderSymbol \equiv a\bracketsInverse{l^{\exponentSymbol}} \bmod m^{\forwardExponentSymbol};
$
the corresponding (\iterateRemainderSymbol, \translationMatrixEntrySymbol)-witness $\left(\exponentSymbol, \addendSymbol\right)$ is such that
$
r \equiv \translationMatrixEntrySymbol\bracketsInverse{-m^{\forwardExponentSymbol}} \bmod l^{\exponentSymbol}.
$

Furthermore, under the assumption that $|\translationValueSymbol| < \max\left(m^{\madicExponentSymbol}, l^{\ladicExponentSymbol}\right)$,  one can show that $\addendSymbol \in \nOpenSet{l^{\exponentSymbol}}$.  We begin by writing
$
\addendSymbol = \frac{l^{\exponentSymbol}\iterateRemainderSymbol - \translationValueSymbol}{m^{\forwardExponentSymbol}}.
$ We can bound the integer $\addendSymbol$ from above with the inequality
\begin{equation}\notag
\addendSymbol < \frac{l^{\ladicExponentSymbol}\left(m^{\madicExponentSymbol} - 1\right) + \max\left(m^{\madicExponentSymbol}, l^{\ladicExponentSymbol}\right)}{m^{\madicExponentSymbol}}  
= 	\begin{cases}
		l^{\ladicExponentSymbol} & \max\left(m^{\madicExponentSymbol}, l^{\ladicExponentSymbol}\right) = l^{\ladicExponentSymbol} \\
		l^{\ladicExponentSymbol} + 1 - \frac{ l^{\ladicExponentSymbol}}{m^{\madicExponentSymbol}} & \max\left(m^{\madicExponentSymbol}, l^{\ladicExponentSymbol}\right) = m^{\madicExponentSymbol}
	\end{cases}  
\end{equation} however, if $\addendSymbol = l^{\exponentSymbol}$, then the equality $m^{\forwardExponentSymbol}\addendSymbol = l^{\exponentSymbol}\iterateRemainderSymbol - \translationValueSymbol$ implies the contradictory equivalence $\translationValueSymbol \equivMod{l^{\exponentSymbol}}0$.    We can bound the integer $\addendSymbol$ from below with the inequality
\begin{equation}\notag
\addendSymbol > \frac{l^{\ladicExponentSymbol} - \max\left(m^{\madicExponentSymbol}, l^{\ladicExponentSymbol}\right)}{m^{\madicExponentSymbol}}  
= 	\begin{cases}
		0 & \max\left(m^{\madicExponentSymbol}, l^{\ladicExponentSymbol}\right) = l^{\ladicExponentSymbol} \\
		-1 + \frac{l^{\ladicExponentSymbol}}{m^{\madicExponentSymbol}} & \max\left(m^{\madicExponentSymbol}, l^{\ladicExponentSymbol}\right) = m^{\madicExponentSymbol};
	\end{cases}  
\end{equation} however, if $\addendSymbol = 0$, then $\translationValueSymbol = l^{\ladicExponentSymbol}\iterateRemainderSymbol \equivMod{l^{\ladicExponentSymbol}} 0$.

We can also show that all such $\left(\iterateRemainderSymbol,\translationMatrixEntrySymbol\right)$-witnesses assume this form. Let $\kappa$ denote the multiplicative order of $l$ in $\integers/ {m^{\forwardExponentSymbol}}\integers$.  If $(\exponentSymbol,\addendSymbol)$ and $(\exponentSymbol',\addendSymbol')$ are such that $e'>e$, and the equalities
$
l^{\exponentSymbol}\iterateRemainderSymbol = m^{\forwardExponentSymbol}\addendSymbol +\translationMatrixEntrySymbol
$ and
$
l^{\exponentSymbol'}\iterateRemainderSymbol = m^{\forwardExponentSymbol}\addendSymbol' + \translationMatrixEntrySymbol
$ hold, then 
$
\left(l^{\exponentSymbol'} - l^{\exponentSymbol}\right)\iterateRemainderSymbol = m^{\forwardExponentSymbol}\left(\addendSymbol' - \addendSymbol\right);
$ consequently, the equivalence
$
l^{\exponentSymbol'}\iterateRemainderSymbol \equiv l^{\exponentSymbol}\iterateRemainderSymbol \bmod m^{\forwardExponentSymbol}
$ holds.   As both $l$ and  $\iterateRemainderSymbol$ are coprime to $m$, we can write
$
l^{\exponentSymbol'-\exponentSymbol} \equiv 1 \bmod m^{\forwardExponentSymbol},
$ and it follows that $\kappa \mid (\exponentSymbol'-\exponentSymbol)$. 

Assume $(\iterateRemainder{v}, \translationMatrixEntry{v}{i})$ is an admissible pair, and let $(\exponentSymbol, \addendSymbol)$ be a $(\iterateRemainder{v}, \translationMatrixEntry{v}{i})$-witness. As $r\equiv i \bmod l$, we will reduce notational clutter by omitting the double subscript when writing $\translationMatrixEntry{v}{i}$, and we will simply write $\translationValue{v}$.

 We identify both the exponent $\exponentSymbol$ and the value $\addendSymbol$ of an $\left(\iterateRemainderSymbol,\translationMatrixEntrySymbol\right)$-witness as a function of the translation value $\translationMatrixEntrySymbol$, the  $m$-residue $\iterateRemainderSymbol$, and a given multiple $h$ of the cyclic order of $l$ in $\integersMod{m^{\forwardExponentSymbol}}$. We will order the $\left(\iterateRemainderSymbol,\translationMatrixEntrySymbol\right)$-witnesses by the exponent element $\exponentSymbol$:  we will say that the witness that includes the $h$-th smallest exponent element is of \emph{height} $h$.

The analysis in this article includes sequences of one or more admissible pairs and witnesses;  we will now define such sequences as follows.

\begin{definition}[Admissible Sequences]  Assume \ref{hypotheses::mlAf}.
For each $w\in \naturals_0$, let $\left(\iterateRemainder{w},\translationValue{v(w)}\right)$ be an admissible pair where $v(w) \equivMod{\tau} w$, and  let $\left(\exponent{w}, \addend{w}\right)$ be the  $\left(\iterateRemainder{w},\translationValue{v(w)}\right)$-witness. 
    We will define the sequence of tuples $\graded{t} = \left(\admissibleFactor{w}\right)_{w\geq 0}$ where
$
 \admissibleFactor{w} = \left(\iterateRemainder{w},\translationValue{v(w)},\height{w}\right)
$ to be an \emph{admissible sequence} of the  system.
\end{definition}

To expedite further reading, we will  establish a second set of hypotheses.

\begin{hypotheses}\label{hypotheses::admissibleParameter} 
Assuming \ref{hypotheses::mlAf},
define $\factor{t}$ to be an admissible sequence of \mlSystem\ of length $\tau$ where $\factor{t} = \left(\admissibleFactor{v}\right)_{v=0}^{\tau-1}$,  and $\admissibleFactor{v} = \left(\admissibleFactorTuple{v}\right)$.   Let $\admissiblePair{v}$ denote the $\left(\iterateRemainder{v}, \translationValue{v}\right)$-witness,  let $\graded{e}$ denote the exponent sequence  $\left(\exponent{v}\right)_{v=0}^{\tau-1}$, and let $\graded{a}$ denote the translation sequence  $\left(\translationValue{v}\right)_{v=0}^{\tau-1}$.   

Furthermore, define $\admissibleFactor{u} = \admissibleFactor{u\bmod \tau}$, and $\admissiblePair{u}= \admissiblePair{u\bmod \tau}$ for $u\in \integers$,  and let $\graded{t}^{\infty} = \left(\admissibleFactor{u}\right)_{u\geq 0}$.
\end{hypotheses}

We will now introduce the \emph{dual-radix mappings} of a given  $(m,l)$-system.  These mappings, which function as the inverse of \graded{l}-adic mappings (when defined),  decide the $\graded{l}$-residue class of its image by conditioning on the $\graded{m}$-residue class of the argument.

\begin{definition}[Dual-Radix Mappings]
Assume \ref{hypotheses::mlAf}.
For each $v \in \nZeroOpenSet{\tau}$, we will define the \emph{dual-radix mapping (of index $v$)} $$\inverseMapFullSymbol{v}: \mathbb{Q} \times\mathbb{N}\times \pm \nZeroOpenSet{m^{\forwardExponent{v}}} \to \mathbb{Q}$$  as follows:  If $\left(\iterateRemainder{v}, \translationValue{v}\right) \in \admissibleResidueSet{v}$,  let $h \in \naturals$, and  let $(\exponent{v}, \addend{v})$ denote the $\left(\iterateRemainder{v}, \translationValue{v}\right)$-witness  of height $h$.  We  define 
 \[
 \inverseMapFull{v}{m^{\forwardExponent{v}}x + \iterateRemainder{v}}{h}{\translationValue{v}} = l^{\exponent{v}}x + \addend{v}.
 \] 
 If   $\left(\iterateRemainder{v}, \translationValue{v}\right) \notin \admissibleResidueSet{v}$, then  define
 \[
\inverseMapFull{v}{m^{\forwardExponent{v}}x + \iterateRemainder{v}}{h}{\translationValue{v}}  =m^{\forwardExponent{v}}x + \iterateRemainder{v}.
 \] 
 
 We define this collection of $\tau$ dual-radix mappings as the \emph{dual-radix system} of the  $(m,l)$-system $\mlSystem$, and we will denote it by $\dualRadixMappings$.

\end{definition}

In the following section, we will turn our attention to infinite, admissible sequences of a given $(m,l)$-system of order $\tau$ of the form $\factor{t}^{\infty}$, where $\factor{t}$ is assumed to be a fixed, admissible sequence of length $\tau$.  \deleted{Following the work of Lagarias in \cite{Lagarias1990}, we will demonstrate that these infinite, periodic admissible sequences correspond to the periodic orbits within the given system, and they give rise to the rational expressions that are the focus of this article.}

\section{The Iterate Values of Periodic Orbits and their Graded Quotients}\label{section::OrbitTerms}

One  open problem in regard to the classic $3x+1$  Problem is the existence of periodic orbits on the domain of positive integers other than the orbit ($1,4,2$).  The work in \cite{BohmSontacchi78} yields the following condition for the existence of such periodic orbits :  an integer is included within an orbit if and only if it admits a rational expression of the form
		\begin{equation}\label{equation::BS}
		\frac{\sum_{0\leq w < \tau}3^w2^{x(\tau-1-w)}}{2^{\sum_{0\leq w < \tau}x(\tau-1-w)}-3^{\tau}};
		\end{equation} here, the term $x(0)=0$,  and the sequence $x(0),\ldots,x(\tau-1)$ is monotonically increasing over $\naturals_0$.   

In this section, we will address a more general question: for a given $(m,l)$-system of order $\tau$, what are the \factor{m}-adic and \factor{l}-adic properties of the iterate values  within a periodic orbit of order $\tau$?  Within the context of the $3x+1$ Problem, such questions have been approached over the extended domains of \rationalUnits{2} (\cite{Lagarias1990}) and $\mathbb{R}$ (\cite{Misiurewicz2005}, \cite{Chamberland1996}).

 First, we will  demonstrate that an iterate $\iterate{}$ within a periodic orbit  assumes a  rational expression analogous to (\ref{equation::BS}) above. Afterwards, we will take a new approach in deciding the \graded{m}-adic and \graded{l}-adic expansions of the iterate: we will use the \factor{l}-adic quotients of each of the $\tau$ iterates to decide the \factor{m}-adic digits of $\iterate{}$, and we will use the \factor{m}-adic quotients of each of the $\tau$ iterates to decide the \factor{l}-adic digits of $\iterate{}$.  
Finally, we will complete the section by defining the relevant terms that arise within our analysis.


\subsection{Iterate Values of a Periodic Admissible Sequence:} 




Under the assumptions in \ref{hypotheses::mlAf} and \ref{hypotheses::admissibleParameter},
we use the information encoded in \graded{t}, to analyze the iterate values of periodic orbits. We begin by establishing the next set of assumptions.
\begin{hypotheses}\label{hypotheses::periodicIterates}  Assume \ref{hypotheses::mlAf} and \ref{hypotheses::admissibleParameter}. Given an admissible sequence \graded{t} of order $\tau$, let $\left(\iterate{v}\right)_{v=0}^{\tau-1}$ denote the sequence of length $\tau$ over  \rationalUnits{m,l}  whose elements  satisfy the equalities
\begin{equation}\label{equation::iterateEquation}
\inverseMapFull{v }{\iterate{v}}{\height{ v}}{\translationValue{ v}} = l^{\exponent{v}}\pfrac{\iterate{v} - \iterateRemainder{v}}{m^{\forwardExponent{v}}} + \addend{v}  = \iterate{  v+1},
\end{equation} where
$
\iterate{v}\equivMod{m^{\forwardExponent{ v}}} \iterateRemainder{{v}}
$ for each $v \in\nZeroOpenSet{\tau}$.  Furthermore, let $\iterate{u} = \iterate{u\bmod \tau}$ for $u\in \integers$.

Also, let $u \in \naturals_0$. We will  define the $u$-th prefix sum $\forwardExponentSum{v,u}$ of $\left[\graded{f} \leftShift  v\right]^{\infty}$ to be
$
\forwardExponentSum{v,u} ~=~\sum_{0\leq y < u}\forwardExponent{ v+y},
$ 
and we will  define the $u$-th prefix sum of $\left[{\graded{e}} \leftShift v\right]^{\infty}$ to be
$
 \leps{v,u}~=~\sum_{0\leq y < u} \lExp{v+y}.
$
We will also define the suffix sums
$
\mess{v,u} = \sum_{0\leq y < u} \mExp{v-1-y},
$ and
$
\less{v,u} = \sum_{0\leq y < u} \lExp{v-1-y}.
$

Define the sequence
$
\numeratorBase{v,u} = \parentheses{m^{\meps{v+u,w}}l^{\less{v,\tau-1-w}}}_{w=0}^{\tau-1},
$  and define
$
\numerator{v} = \numeratorBase{v,0} \cdot \brackets{\graded{a}\leftShift v}.
$
Finally, define $\denominator{v} = l^{\less{v,\tau}} - m^{\meps{v,\tau}}$.
\end{hypotheses}


 We will demonstrate the sequence of iterates $\left(\iterate{v}\right)_{v=0}^{\tau-1}$ exists, and its $v$-th element satisfies the equality
		\begin{equation}\label{equation::IterateRationalExpression}
		\iterate{ v} = \frac{\numerator{v}}{\denominator{v}}
		\end{equation} for each $v\in\nZeroOpenSet{\tau}$; this identity was established in the $3x+1$ context in \cite{BohmSontacchi78}.

Before doing so, we will establish some elementary properties of  the sums. 

\begin{identities}\label{identities::ExponentSums} For $v\in \nZeroOpenSet{\tau}$ and $u\in \integers$, the following identities hold for the prefix and suffix sums of $\brackets{\graded{f}\leftShift v}^{\infty}$:
	\begin{enumerate}[i.]
		\item\label{ExponentSumsShift}  $\meps{v,\tau} = \meps{uv,\tau} = \mess{uv,\tau}$;
		\item\label{ExponentSumsReverseSplit} for $w\in \nZeroSet{\tau}$, we have $\meps{v+w,\tau-w}= \mess{v,\tau-w}$; 
		\item\label{ExponentSumsSplit} for $w\in \nZeroSet{\tau}$, we have $\meps{v,\tau}  = \meps{v,w} +\meps{v+w,\tau-w}$;
		\item\label{PrefixExponentSumsPerturb}	$\meps{v,u} + \mExp{v+u} = \meps{v,u+1} = \mExp{v} + \meps{v+1,u}$;
		\item\label{SuffixExponentSumsPerturb}	 $\mess{v,u} + \mExp{v-1-u} = \mess{v,u+1} = \mExp{v-1} + \mess{v+1,u}.$	
	\end{enumerate}
	
	Replacing $\madicExponentSumSymbol$ with  $\ladicExponentSumSymbol$ and  $\madicExponentSymbol$ with $\ladicExponentSymbol$ yields analogous identities  for the prefix and suffix sums of $\brackets{\graded{e}\leftShift v}^{\infty}$.
\end{identities}

Our approach to demonstrating (\ref{equation::IterateRationalExpression}) will be similar to the approach presented in \cite{BohmSontacchi78};  we will appeal to the $\graded{l}$-adic maps of the system.  For all $v\in \nZeroOpenSet{\tau}$, the condition in  (\ref{equation::iterateEquation}) is equivalent to the condition 
$$\forwardMap{\iterate{  v+1}}{v} = \frac{m^{\forwardExponent{ v}}\iterate{  v+1} + \translationValue{v}}{l^{\exponent{v}}}  = \iterate{ v}.$$

 First, we will show by induction that the iterate \iterate{ v} satisfies the equation
\begin{equation}\label{equation::iterateInduction}
\iterate{ v} = \pfrac{m^{\meps{v,u}} }{l^{\leps{v,u}}}\iterate{v+u} + \sum_{0 \leq w < u} \pfrac{m^{\meps{v,w} }}{l^{\leps{v,w+1}}}\translationValue{v+w}
\end{equation}  for $u\in \naturals_0$.    The claim readily follows for $u=0$. Assuming the claim,  we will apply the substitution
\[
\iterate{v+u} =  \frac{m^{\mExp{v+u}}\iterate{v+u+1} + \translationValue{v+u}}{l^{\lExp{v+u}}}  
\] into $\left(\ref{equation::iterateInduction}\right)$.

By induction, along with the identities (\ref{identities::ExponentSums}.\ref{PrefixExponentSumsPerturb}) and (\ref{identities::ExponentSums}.\ref{SuffixExponentSumsPerturb}), we arrive at the equality
\[
\iterate{ v} = \pfrac{m^{\meps{v, u+1}} }{l^{\leps{v, u+1}}}\iterate{v+u+1} + \sum_{0 \leq w < u+1} \pfrac{m^{\meps{v,w} }}{l^{\leps{v,w+1}}}\translationValue{v+w}
\]
as required.

The condition   $\iterate{v+\tau} = \iterate{ v}$  yields the equality
\[
\iterate{ v} =  \pfrac{m^{\meps{v, \tau}}}{l^{\leps{v, \tau}}}\iterate{ v} + \sum_{0\leq w < \tau} \pfrac{m^{\meps{v,w}}}{l^{\leps{v,w+1}}}\translationValue{v+w}
\] when $u=\tau-1$. Solving for the iterate $\iterate{ v} $,  by way of the identities
 in (\ref{identities::ExponentSums}.\ref{ExponentSumsSplit}),
 yields the equality
\[
\iterate{ v} =  \frac{\sum_{0\leq w < \tau} m^{\meps{v,w}}l^{\less{v,\tau-1-w}}\translationValue{v+w}   }{l^{\less{v,\tau}} - m^{\meps{v, \tau}}}.
\]  

As $m$ and $l$ are coprime, the difference $l^{\less{v,\tau}} - m^{\meps{v,\tau}}$ is coprime to both $m$ and $l$;  consequently,  this rational formulation of \iterate{v} guarantees that it admits expansions within  $Z_{\graded{m},\graded{f}\leftShift v}$ and $Z_{\graded{l},\graded{e}^R\rightShift v}$.   We will turn our focus to the $\factor{m}$-adic  and $\factor{l}$-adic digits  of  \iterate{v} (over the respective rings); these digits will be determined inductively using  the dual-radix mappings of the given $(m,l)$-system. 

\subsection{The \factor{m}-adic and \factor{l}-adic Digits of the Graded Quotients}

We continue our analysis of \iterate{v} from Assumptions \ref{hypotheses::periodicIterates}. We will proceed with deciding the \factor{m}-adic digits of \iterate{v} by considering the dual-radix system of the given $(m,l)$-system.

\subsubsection{The \factor{m}-adic Digits of \iterate{ v} over $Z_{\graded{m},\graded{f}\leftShift v}$}
Let $v\in \nZeroOpenSet{\tau}$. 
For each $u\in \naturals_0$, we will define $\madicShift{v,u}$ to be the $u$-th $\factor{m}$-adic quotient of $\iterate{v}$ over $Z_{\graded{m},\graded{f}\leftShift v}$: we will assume that $\madicShift{v,u}$ satisfies the equality
\begin{equation}\label{equation::madicShiftExpand}
\madicShift{v,u}= m^{\mExp{v+u}}\madicShift{v,u+1} + \madicDigit{v,u} 
\end{equation} where $\madicDigit{v,u} \in \nZeroOpenSet{m^{\mExp{v+u}}}$.   

We will inductively determine the $\factor{m}$-adic digits $\left(\madicDigit{v,u}\right)_{u\geq 0}$ of \iterate{ v} over $Z_{\graded{m},\graded{f}\leftShift v}$ by way of the identity
\[
\iterate{v} = m^{\mExp{v}}\madicShift{v,1} + \iterateRemainder{ v} 
\xrightarrow{\inverseMapSymbol} l^{\exponent{ v  }}\madicShift{v,1} + \addend{ v} \notag \\
= \iterate{v+1}.
\]   In doing so, we will express the quotient \madicShift{v+u,1} with the equality
\[
		\madicShift{v+u,1} = l^{\leps{v,u}}\madicShift{v,u+1} +\ladicPrefixAddend{v,u}
\] where $\ladicPrefixAddend{v,u} \in \integers$ (we  define $\ladicPrefixAddend{v,0} = 0$). We will also show that both $\madicDigit{v,u}$ and $\ladicPrefixAddend{v,u}$ are functions of the $u+1$ terms $\left(\admissibleFactor{v+w}\right)_{w=0}^{u}$ of $\graded{t}^{\infty}$.

When $u = 0$, the inductive claim holds as $\iterateRemainder{v} = \madicDigit{v,0} = \madicDigit{v,0}\left(\admissibleFactor{v}\right)$, and $\ladicPrefixAddend{v,0} = 0$.

For $u=1$, we begin by writing $\iterate{ v} =m^{\mExp{v}}\madicShift{v,1} + \iterateRemainder{ v}$, and we  will write
\[
\iterate{v+1} = l^{\lExp{v}}\madicShift{v,1} + \addend{ v} = m^{\mExp{v+1}}\madicShift{v+1,1} +\iterateRemainder{v+1}.
\]  We will express the first \factor{m}-adic quotient  \madicShift{v,1} of $\iterate{ v}$ in terms of its \factor{m}-adic quotient (over $Z_{\graded{m},\graded{f}\leftShift v+1}$) and its first \factor{m}-adic digit: we will write  $\madicShift{v,1} = m^{\mExp{v+1}}\madicShift{v,2} + \madicDigit{v,1}$ where $\madicDigit{v,1} \in \nZeroOpenSet{m^{\mExp{v+1}}}$.   Thus, as $\lExp{v} = \leps{v,1}$, we will write
\[
l^{\lExp{v}}\left(m^{\mExp{v+1}}\madicShift{v,2} +\madicDigit{v,1} \right)+  \addend{ v}  = m^{\mExp{v+1}}\left(l^{\leps{v,1}}\madicShift{v,2} + \ladicPrefixAddend{v,1} \right) + \iterateRemainder{v+1} = \iterate{v+1},
\]
where the term 
\[
 \ladicPrefixAddend{v,1}  =  \frac{l^{\leps{v,1}}\madicDigit{v,1}- \iterateRemainder{v+1} +\addend{ v}}{m^{\mExp{v+1}}}.
\] 

We can now express the value $\madicShift{v+1,1}$ with the equality
\begin{equation}\label{equation::prefixAddendClosure}
\madicShift{v+1,1} = l^{\leps{v,1}}\madicShift{v,2} + \ladicPrefixAddend{v,1}.
\end{equation} Under the assumption that both $\iterate{v}$ and $\iterate{v+1}$ are elements of  $\rationalUnits{m,l}$, it follows that the \graded{m}-adic quotients $\madicShift{v,2}$ and $\madicShift{v+1,1}$   are  elements of $\rationalUnits{m,l}$.  Furthermore, from the equality in (\ref{equation::prefixAddendClosure}), we conclude that the term $ \ladicPrefixAddend{v,1} = \madicShift{v+1,1} -  l^{\leps{v,1}}\madicShift{v,2} $ is also an element of $\rationalUnits{m,l}$ (as per the ring closure properties); thus, we require that 
$
m^{\mExp{v+1}} \pdiv \left[l^{\leps{v,1}}\madicDigit{v,1}- \iterateRemainder{v+1} + \addend{ v}\right].
$
 This  requirement determines the unique \\$m^{\mExp{v+1}}$-residue 
\[
\madicDigit{v,1} = \madicDigit{v,1} \left(\admissibleFactor{v}, \admissibleFactor{ v+1}\right) \equivMod{m^{\mExp{v+1}}} \left[l^{\lExp{v}}\right]^{-1}\left(\iterateRemainder{v+1} - \addend{ v}\right).
\]   Consequently, we can now write $ \ladicPrefixAddend{v,1} = \ladicPrefixAddend{v,1}\left(\admissibleFactor{v}, \admissibleFactor{ v+1}\right)$.

One can proceed inductively on $u$ and,  by appealing  (\ref{equation::madicShiftExpand}) for each $u$, express the iterate \iterate{v+u} as
		\begin{equation}\label{equation::backwardIterate}
		       \iterate{v+u} = m^{\mExp{v+u}}\madicShift{v+u,1} + \iterateRemainder{v+u}
		\end{equation} 
		where 
		\begin{equation}\label{equation::backwardQuotient}
		\madicShift{v+u,1} = l^{\leps{v,u}}\madicShift{v,u+1} +\ladicPrefixAddend{v,u},
		\end{equation} and
		\[
		\ladicPrefixAddend{v,u} =  \frac{l^{\leps{v,u}}\madicDigit{v,u} + l^{\lExp{v+u-1}}\ladicPrefixAddend{v, u-1} - \iterateRemainder{v+u} +\addend{v+u-1}
}{m^{\mExp{v+u}}}.
		\]  Similar to the above, the condition that both $\madicShift{v+u,1}$ and $\madicShift{v,u+1}$ are elements of $\rationalUnits{m,l}$ requires that the condition $$m^{\mExp{v+u}} \pdiv  \left[l^{\leps{v,u}}\madicDigit{v,u} + l^{\lExp{v+u-1}}\ladicPrefixAddend{v,u-1} - \iterateRemainder{v+u} +\addend{v+u-1}\right]$$   is met;  this requirement determines the unique $m^{\mExp{v+u}}$-residue 
		\begin{equation}\label{equation::quotientRemainderEquivalence}
		\madicDigit{v,u} \equivMod{m^{\mExp{v+u}}} \left[l^{\leps{v,u}}\right]^{-1}\left[ -  l^{\lExp{v+u-1}}\ladicPrefixAddend{v,u-1} + \iterateRemainder{v+u} - \addend{v+u-1}\right].
		\end{equation}   Furthermore, this equivalence, along with the inductive assumption, allows us to conclude that 
		$
		\madicDigit{v,u} = \madicDigit{v,u} \left(\admissibleFactor{v},\ldots, \admissibleFactor{v+u}\right),
		$ and 
		$
		  \ladicPrefixAddend{v,u}  = \ladicPrefixAddend{v,u}\left(\admissibleFactor{v},\ldots, \admissibleFactor{v+u}\right).
		$

\subsubsection{The \factor{l}-adic Digits of \iterate{ v} over $\padics{\graded{l},\graded{e}^R\rightShift v}$}
  Let $v\in \nZeroOpenSet{\tau}$.  
For each $u\in \naturals_0$, we will define $\ladicShift{v,u}$ to be the $u$-th $\factor{l}$-adic quotient of $\iterate{v}$ over $\padics{\graded{l},\graded{e}^R\rightShift v}$:  we will assume that $\ladicShift{v,u}$ satisfies the equality
$
\ladicShift{v,u} = l^{\lExp{v-1-u}}\ladicShift{v,u+1} + \ladicDigit{v,u}
$ where $\ladicDigit{v,u}\in \nZeroOpenSet{l^{\lExp{v-1-u}}}$. 

		 Analogous to above, one can inductively  determine the $\factor{l}$-adic digits of \iterate{ v} over $\padics{\graded{l},\graded{e}^R\rightShift v}$  by way of the identity
		 \[
		 \iterate{ v} =  l^{\lExp{v-1}}\ladicShift{v,1} + \addend{v-1} 
		 \xrightarrow{\inverseMapSymbol^{-1}} m^{\mExp{v-1}}\ladicShift{v,1}  + \iterateRemainder{v-1} 
		 = \iterate{v-1} 
		\]  for each $v\in \nZeroOpenSet{\tau}$.
Induction allows one to express the quotient \ladicShift{v-u,1} with the equality
		$
		\ladicShift{v-u,1} = m^{\mess{v,u}}\ladicShift{v,u+1} +\madicPrefixAddend{v,u},
		$ where 
				\[
\madicPrefixAddend{v,u} = \frac{m^{\mess{v,u}   }\ladicDigit{v,u}+ m^{\mExp{v-u}}\madicPrefixAddend{v,u-1}+ \iterateRemainder{v-u}  -\addend{v-1-u}  }{l^{\lExp{v-u-1}}}
\] 
(we  define $\madicPrefixAddend{v,0} = 0$), and
\begin{equation}\label{equation::diagonalDifferenceEquivalence}
\ladicDigit{v,u}\equivMod{l^{\lExp{v-1-u}}} \left[-m^{\mess{v,u}  }\right]^{-1}\left[ m^{\mExp{v-u}}\madicPrefixAddend{v,u-1} + \iterateRemainder{v-u}  -\addend{v-u-1}   \right].
\end{equation}  Furthermore, we have that 
		$
		\ladicDigit{v,u}= \ladicDigit{v,u}\left(\admissibleFactor{v-1},\ldots, \admissibleFactor{ v-u-1}\right),
		$ and \\
		$
		  \madicPrefixAddend{v,u}  = \madicPrefixAddend{v,u}\left( \admissibleFactor{v-1},\ldots, \admissibleFactor{v-u-1}\right).
		$

		We have shown that, for all $u\geq 0$, the $u$-th \factor{m}-adic and \factor{l}-adic digits of the iterate $\iterate{v}$ within a periodic orbit are uniquely decided by the admissible factor $\graded{t} \leftShiftOperator{v}$; we will summarize  the terms we have encountered thus far in the following definition.

\begin{definition}\label{definition::mlSystemsTerms}
Assume \ref{hypotheses::periodicIterates}.
For $u\in \naturals_0$, 
		\begin{enumerate}
		        
		        \item let \madicShift{v,u} and \madicDigit{v,u} denote the $u$-th \factor{m}-adic quotient and \factor{m}-adic digit, respectively, of $\iterate{ v}$ over $\padics{\graded{m},\graded{f}\leftShift v}$. We can write	        
		$
		\madicShift{v+u,1}~=~l^{\leps{v,u}}\madicShift{v,u+1} + \ladicPrefixAddend{v,u}  ,
		$ where  the term $ \ladicPrefixAddend{v,u} $, the $u$-th \factor{l}-adic \emph{prefix addend} of $\iterate{ v}$, is to be integer-valued, and it is defined recursively as
		        \begin{equation}
		                \ladicPrefixAddend{v,u} = 
		                \begin{cases}
		                        0 & u = 0,\\
		                      \frac{l^{\leps{v,u}}\madicDigit{v,u} + l^{\lExp{v+u-1}} \ladicPrefixAddend{v,u-1}  - \iterateRemainder{v+u} +\addend{v+u-1}
}{m^{\mExp{v+ u }}}  & u > 0. \notag
		                \end{cases}
		        \end{equation}
		\item let \ladicShift{v,u} and \ladicDigit{v,u} denote the $u$-th \factor{l}-adic quotient and \factor{l}-adic digit, respectively,  of $\iterate{ v}$ over $\padics{\graded{l},\graded{e}^R\rightShift v}$. We can write
		$
	\ladicShift{v-u,1}~=~m^{\mess{v,u}}\ladicShift{v,u+1} +\madicPrefixAddend{v,u},
		$ where  the term $\madicPrefixAddend{v,u}$, the $u$-th \factor{m}-adic \emph{prefix addend} of $\iterate{ v}$, is to be integer-valued, and it is defined recursively as
		        \begin{equation}
		              \madicPrefixAddend{v,u} = 
		                \begin{cases}
		                        0 & u = 0,\\
		                         \frac{m^{\mess{v,u}   }\ladicDigit{v,u}+ m^{\mExp{v-u}}\madicPrefixAddend{v,u-1}+ \iterateRemainder{v-u}  -\addend{v-u-1}  }{l^{\lExp{v-1-u}}} & u > 0. \notag
		                \end{cases}
		        \end{equation}		        
		        \end{enumerate}

		\end{definition}

		We have shown that the iterate $\iterate{v}$ can be expressed as
		\[
		\iterate{ v} =   \frac{\numeratorBase{v,0}\cdot\brackets{\graded{a}\leftShift v}}{\denominator{v}}
		\] for all $v\in \nZeroOpenSet{\tau}$;  we will call this the \emph{dual-radix rational form} of the iterate \iterate{v}\replaced{.}{, and we will examine the \factor{m}-adic and \factor{l}-adic properties of such expressions.} In the following section, we will demonstrate that the \factor{m}-adic and \factor{l}-adic quotients of these  iterate values  assume a similar rational form:  for  $u\geq 0$, we can express the $u$-th \factor{m}-adic quotient $\madicShift{v,u}$ of $\iterate{v}$ (over $Z_{\graded{m},\graded{f}\leftShift v}$) as
		\[
		\madicShift{v,u} = \frac{\numeratorBase{v,u}\cdot\brackets{\graded{c}_u\leftShift v}}{\denominator{v}}
		\] where
		\begin{equation}
		\digitDifference{v+w,u} =
			\begin{cases}
				\translationValue{v+w} & u=0 \\
				\madicDigit{v+w+1,u-1} -\ladicDigit{v+w+u,u-1} & u>0.
			\end{cases}\notag
		\end{equation}
%
		 Furthermore, we will connect the \factor{m}-adic and \factor{l}-adic quotients of the iterates via the equality
		$
		\ladicShift{v+u,u} = \madicShift{v,u}.
		$  
		We will also demonstrate how to  compute efficiently the \graded{m}-adic and \graded{l}-adic digits of index $u$ of all $\tau$ iterates from the differences $\left(\digitDifference{v,u}\right)_{v\in \nZeroOpenSet{\tau}}$.  
		

\section{Results}\label{section::MainResults}
The first theorem of this article pertains to the arithmetic differences of the prefix addends of the form
$
\ladicPrefixAddend{v,u}  - \ladicPrefixAddend{v+1,u-1} 
$ and
$
\madicPrefixAddend{v,u}  - \madicPrefixAddend{v-1, u-1}.
$ 
 In the previous section, we have demonstrated the role of the \graded{l}-adic prefix addend $\ladicPrefixAddend{v,u}$  in deciding the $\graded{m}$-adic digit $\madicDigit{v,u}$ of the iterate $\iterate{v}$.  We will demonstrate that the difference $\ladicPrefixAddend{v,u}  - \ladicPrefixAddend{v+1,u-1} $ yields the $\factor{l}$-adic digit $\ladicDigit{ v+u+1, u}$ of the iterate $\iterate{v+u+1}$. Similarly, we have shown that the \graded{m}-adic prefix addend $\madicPrefixAddend{v,u}$  decides the $\graded{l}$-adic digit $\ladicDigit{v,u}$; we will now show that the  difference $\madicPrefixAddend{v,u}  - \madicPrefixAddend{v-1, u-1}$ yields the $\factor{m}$-adic digit $\madicDigit{v-u-1,u}$ of the iterate $\iterate{v-u-1}$. One consequence of this result is  that the $\factor{m}$-adic image of the \graded{m}-adic quotient $\madicShift{v,1}$ $\left(\textmd{over } \padics{\graded{m},\graded{f}\leftShift v+1}\right)$ is 
$
 \left(\madicPrefixAddend{v+u+1,u}\right)_{u\geq 1}, 
$ and the $\factor{l}$-adic image  of the \graded{l}-adic quotient $\ladicShift{v,1}$ $\left(\textmd{over } \padics{\graded{l},\graded{e}^R\rightShift v-1}\right)$ is  
$
 \left(\ladicPrefixAddend{v-u-1,u}\right)_{u\geq 1}.
$

As we will see later in this section, this {\sl dual-radix property} of the prefix addends--- simultaneously deciding both an \graded{m}-adic and \graded{l}-adic digit of an iterate---gives rise to a  recurrence for computing  the digits (in both bases) of all the iterates within a periodic orbit simultaneously.

\begin{theorem}[\factor{l}-adic Diagonal Differences]\label{theorem:diagonalDifferences}  
Assume \ref{hypotheses::periodicIterates} and \ref{definition::mlSystemsTerms}. 

 For each $v\in \nZeroOpenSet{\tau}$ and $u\in \naturals$, the
prefix addend difference
        $
                \ladicPrefixAddend{v,u}  - \ladicPrefixAddend{v+1,u-1} 
        $ equals
      	$
          l^{\less{v+u,u-1}}\ladicDigit{v+u+1,u};
        $  furthermore, the \factor{l}-adic digit $\ladicDigit{v+u+1,u}$ admits the recurrence  
		\begin{equation}\label{equation::digitRecurrence}        
\ladicDigit{v+u+1,u} =
\frac{l^{\lExp{v}}\madicDigit{v,u}  - \madicDigit{v+1,u-1}+
\ladicDigit{v+u,u-1} }{m^{\mExp{v+u}}}.
        \end{equation}
\end{theorem}

\begin{proof}

	Assume the hypotheses and notation within the theorem statement. We will demonstrate the claim by induction on $u$.
	In order to do so, we will establish the identity
	\begin{equation}\label{equation::LAdicPrefixExpansion}
	\ladicPrefixAddend{v,u}  = \sum_{0\leq w < u} l^{ \less{v+u,w}} \ladicDigit{v+u+1,w+1}
	\end{equation} for all $v\in \nZeroOpenSet{\tau}$ and $u\in \naturals$.

Firstly, the difference 
\begin{equation}\label{equation::lemmaBaseCase}
\ladicPrefixAddend{v,1} - \ladicPrefixAddend{v+1,0} = \ladicPrefixAddend{v,1} = \frac{l^{\lExp{v}}\madicDigit{v,1} -\iterateRemainder{v+1}+\addend{v} }{m^{\mExp{v+1}}};
\end{equation}  
we will show that this difference of integers is an  element of $\nZeroOpenSet{l^{\lExp{v}}}$: one  can bound this difference from below with the inequalities
\[
\ladicPrefixAddend{v,1} > -\frac{m^{\mExp{v+1}}-1}{m^{\mExp{v+1}}} > -1;
\] furthermore, as the admissible addend $\addend{v} \in \nZeroOpenSet{l^{\lExp{v}}}$, we can bound the difference in (\ref{equation::lemmaBaseCase})  from above with the inequalities
\[
\ladicPrefixAddend{v,1} < \frac{l^{\lExp{v}}\left(m^{\mExp{v+1}} - 1\right) + l^{\lExp{v}} - 1}{m^{\mExp{v+1}}} <  l^{\lExp{v}}.
\]
As $
\iterate{ v+1} =  l^{\lExp{v}}\ladicShift{v+1,1} + \addend{ v} =  \inverseMapFull{v}{\iterate{v}}{\height{v}}{\translationValue{v}} = l^{\lExp{v}}\madicShift{v,1} + \addend{ v},
$ we have the equality 
$
\ladicShift{v+1,1} = \madicShift{v,1}
$  for each $v\in \nZeroOpenSet{\tau}$
As we can write
$
  \ladicShift{v+2,1} =  \madicShift{v+1,1} = l^{\less{v+1,1}}\madicShift{v,2}  + \ladicPrefixAddend{v,1}
$ as per (\ref{equation::backwardQuotient}),  the fact that $\ladicPrefixAddend{v,1} \in \nZeroOpenSet{l^{\lExp{v}}}$
yields the equality
$
\ladicPrefixAddend{v,1} = \ladicDigit{v+2,1}.
$ 

We will now proceed inductively on $u$.  
We will appeal to (\ref{identities::ExponentSums}.\ref{PrefixExponentSumsPerturb}) and express the difference  $\ladicPrefixAddend{v,u+1}  - \ladicPrefixAddend{v+1,u}$ as
\begin{equation}\label{equation:prefixAddendDifference}
\frac{l^{\leps{v+1,u}}\left[l^{\lExp{v}}\madicDigit{v,u+1} - \madicDigit{v+1,u}\right] + l^{\lExp{v+u}}\left(\ladicPrefixAddend{v,u}  - \ladicPrefixAddend{v+1,u-1} \right) }{m^{\mExp{v+u+1}}}.
\end{equation}  By induction, the difference 
$
\ladicPrefixAddend{v,u}  - \ladicPrefixAddend{v+1,u-1}  =  l^{\less{v+u,u-1}}\ladicDigit{v+u+1,u}.
$ 
 We can complete the inductive argument by writing $\less{v+u+1,u} = \leps{v+1,u} = \lExp{v+u} +\less{v+u,u-1}$, and rewriting (\ref{equation:prefixAddendDifference})
 as
\[
l^{\less{v+u+1,u}}\bfrac{l^{\lExp{v}}\madicDigit{v,u+1} - \madicDigit{v+1,u} +\ladicDigit{v+u+1,u}}{m^{\mExp{v+u+1}}}.
\]  Similar to the reasoning above, the value  of 
\[
\frac{l^{\lExp{v}}\madicDigit{v,u+1} - \madicDigit{v+1,u} +\ladicDigit{v+u+1,u}}{m^{\mExp{v+u+1}}}
\] is an element of \nZeroOpenSet{l^{\lExp{v}}}.  We denote the value of this expression by $\beta$, and we will write
\begin{alignat}{2}
\ladicShift{v+u+2,1} & = \madicShift{v+u+1,1}\notag \\ 
&=  l^{\leps{v,u+1}}\madicShift{v,u+2} + \ladicPrefixAddend{v,u+1} \notag \\
&= l^{\leps{v,u+1}}\madicShift{v,u+2}  + l^{\less{v+u+1,u}}\beta+ \ladicPrefixAddend{v+1,u}.\label{equation::ladicInduction}
\end{alignat} By induction, the prefix addend  $\ladicPrefixAddend{v+1,u}$ is an $l^{\less{v+u+1,u}}$-residue; thus, $\ladicDigit{v+2+u,u+1}~=~\beta$, and
$$
\ladicPrefixAddend{v,u+1} = \sum_{0\leq w < u+1} l^{ \less{v+u+1,w}} \ladicDigit{v+u+2,w+1}
$$

\end{proof}

\begin{corollary}[\factor{m}-adic Diagonal Differences] \label{corollary::madicDiagonalDifferences}

 For $v\in \nZeroOpenSet{\tau}$ and $u\in\naturals$, the difference of $\factor{m}$-adic prefix addends $\madicPrefixAddend{v,u} - \madicPrefixAddend{v-1,u-1}$ equals $m^{\forwardExponentSum{v-u,u-1}} \madicDigit{v-u-1,u}, 
$ where
\[
\madicDigit{v-u-1,u} = \frac{m^{\mExp{v-1}} \ladicDigit{v,u}  +\madicDigit{v-u,u-1} - \ladicDigit{v-1,u-1} }{l^{\lExp{v-u-1}}}.
\]
\end{corollary}

\begin{proof} Assume the notation and hypotheses in the statement of the corollary. 

The proof, similar to the above, is by induction on $u$: in this case, we  express the difference $\madicPrefixAddend{v,u+1} -\madicPrefixAddend{v-1,u}$ as
\[
 m^{\meps{v-u-1,u}}\left[  \frac{m^{\mExp{v-1}}\ladicDigit{v,u+1} + \madicDigit{v-u-1,u} - \ladicDigit{v-1,u}    }{l^{\lExp{v-u-2}}} \right].
\] As a result of Theorem \ref{theorem:diagonalDifferences}, we can express this difference as
$
m^{\meps{v-u-1,u}} \madicDigit{v-u-2,u+1}.
$
\end{proof}

The proof of the theorem also allows us to equate the \factor{m}-adic and \factor{l}-adic quotients of the orbit elements.

\begin{corollary} 
\label{corollaryQuotientIdentity} For  $v\in \nZeroOpenSet{\tau}$ and $u\in \naturals_0$, the equality
        $
		\ladicShift{v+u,u} = \madicShift{v,u}
	$ holds.
\end{corollary}

\begin{proof}  Fix  $v\in \nZeroOpenSet{\tau}$. The claim follows from the definitions when $u=0$.
For $u \in \naturals$, we can write
 $
  \ladicShift{v+u,1}  = l^{\leps{v,u-1}}\ladicShift{v+u,u} +\sum_{0\leq w < u-1} l^{ \less{v+u-1,w}} \ladicDigit{v+u,w+1}$; as $\ladicShift{v+u,1} = \madicShift{v+u-1,1} $, we can also write  $  \ladicShift{v+u,1}= l^{\leps{v,u-1}}\madicShift{v,u} +\ladicPrefixAddend{v,u-1} 
 $ as per (\ref{equation::backwardQuotient}).  As $\ladicPrefixAddend{v,u-1}  = \sum_{0\leq w < u-1} l^{ \less{v+u-1,w}} \ladicDigit{v+u,w+1}$, the result follows.
\end{proof}

The recurrence in (\ref{equation::digitRecurrence}) gives rise to a novel recurrence for computing the canonical expansions of $\iterate{v}$, over  $\padics{\graded{m},\graded{f}\leftShift v}$ and $ \padics{\graded{l},\graded{e}^R\rightShift v}$, for each $v\in \nZeroOpenSet{\tau}$,  simultaneously.

\begin{theorem}\label{theorem::quotientDigitRecurrence}  
Assume \ref{hypotheses::periodicIterates} and  \ref{definition::mlSystemsTerms}. 

For $v\in \nZeroOpenSet{\tau}$, the digits $\madicDigit{v,u}$ and $\ladicDigit{v,u}$ satisfy the following recurrence: for $u\in \naturals$,   we express the $\factor{m}$-adic digit \madicDigit{v,u} by
\[
	\madicDigit{v,u} \equivMod{m^{\mExp{v+u}}}
\left[l^{\lExp{v}}\right]^{-1}\left[ \madicDigit{v+1,u-1} - \ladicDigit{v+u,u-1}\right],
	\]  and  the $\factor{l}$-adic digit \ladicDigit{v,u} by
	\[
		\ladicDigit{v,u} \equivMod{l^{\lExp{v-u-1}}} \left[-m^{\mExp{v-1}}\right]^{-1}\left[ \madicDigit{v-u,u-1} - \ladicDigit{v-1,u-1}\right],
	\]  
	along with the initial conditions $\madicDigit{v,0} = \iterateRemainder{v}$ and
$\ladicDigit{v,0}= \addend{v-1}$. 
\end{theorem}

\begin{proof} Assume the hypotheses and notation within the theorem statement.
 
From Theorem \ref{theorem:diagonalDifferences},  the $\factor{m}$-adic digit  
 \[
 \madicDigit{v,u} = \frac{m^{\mExp{v+u}}\ladicDigit{v+u+1,u} +\madicDigit{v+1,u-1} - \ladicDigit{v+u,u-1} }{l^{\lExp{v}}},
 \] and the $\factor{l}$-adic digit
  \[
 \ladicDigit{v,u} = \frac{l^{\lExp{v-u-1}}\madicDigit{v-u-1,u} - \madicDigit{v-u,u-1} +\ladicDigit{v-1,u-1}}{m^{\mExp{v-1}}}
 \] for $u\in\naturals$.
  If we rewrite the first equation as 
 \[
  l^{\lExp{v}}\madicDigit{v,u} = m^{\mExp{v+u}} \ladicDigit{v+u+1,u} +
\madicDigit{v+1,u-1} -\ladicDigit{v+u,u-1}, 
 \] then we arrive at the equivalence 
$
 \madicDigit{v,u} \equivMod{m^{\mExp{v+u}}}
\left[l^{\lExp{v}}\right]^{-1}\left[ \madicDigit{v+1,u-1} - \ladicDigit{v+u,u-1}\right].
$
 
  If we rewrite the second equation as 
 \[
m^{\mExp{v-1}} \ladicDigit{v,u}    = l^{\lExp{v-u-1}}\madicDigit{v-u-1,u}-
\madicDigit{v-u,u-1} +\ladicDigit{v-1,u-1}, 
 \]  then we arrive at the equivalence 
 $
		\ladicDigit{v,u} \equivMod{l^{\lExp{v-u-1}}} \left[-m^{\mExp{v-1}}\right]^{-1}\left[ \madicDigit{v-u,u-1} - \ladicDigit{v-1,u-1}\right].
	$	
  \end{proof}


For the remainder of this article, the differences of \graded{m}-adic and \graded{l}-adic digits that arose in the previous theorem are encapsulated in the following defintion.

\begin{definition}\label{definition::digitDifference} Assume  \ref{hypotheses::periodicIterates} and  \ref{definition::mlSystemsTerms}.  

For $v\in \nZeroOpenSet{\tau}$ and $u\in \naturals_0$, let $\digitDifference{v,u}$  denote the digit difference
\begin{equation}
		\digitDifference{v,u} =
			\begin{cases}
				\translationValue{v} & u=0 \\
				\madicDigit{v+1,u-1} -\ladicDigit{v+u,u-1} & u>0.
			\end{cases}\notag
		\end{equation}  
		We will also write $\crossSection{v,u} = \crossSection{u}\leftShift v=  \left(\digitDifference{v+w,u}\right)_{w=0}^{\tau-1}$.
\end{definition}
 Theorem \ref{theorem::quotientDigitRecurrence} uses these digit differences  to  compute recursively both the $\factor{m}$-adic and $\factor{l}$-adic digits for each of the $\tau$ iterates within the periodic orbit;  we will define the {\sl quotient cylinder} of the orbit to be the doubly-indexed sequence $\left(\digitDifference{v,u}\right)$ where $v\in \nZeroOpenSet{\tau}$ and $ u\in \naturals_0$. Theorem \ref{theorem::mixedRadixQAdicQuotientShift} shows that the $u$-th cross-section of the quotient cylinder  $ \crossSection{u} = \left(\digitDifference{w,u}\right)_{\substack{w\in \nZeroOpenSet{\tau}}}$ yields a method for expressing the \factor{m}-adic and \factor{l}-adic quotients of index $u$ of the iterate values as dual-radix rationals.

\begin{theorem}\label{theorem::mixedRadixQAdicQuotientShift} 
Assume \ref{hypotheses::periodicIterates},   \ref{definition::mlSystemsTerms}, and \ref{definition::digitDifference}. 

%
Then, for each $v\in \nZeroOpenSet{\tau}$, the equality 
	$$
	\madicShift{v,u}= \frac{\numeratorBase{v,u}\cdot \brackets{\crossSection{u}\leftShift v} }{\denominator{v}} 
	$$ holds.
\end{theorem}

\begin{proof} 

Assume the hypotheses and notation in the theorem statement.  We will prove the claim by induction on $u$. 

We will demonstrate the equality
$
 \denominator{v}\left(\madicShift{v,u} - \madicDigit{v,u}\right) = m^{\mExp{v+u}}\numeratorBase{v,u+1}\cdot \brackets{\crossSection{u+1}\leftShift v}
$ for all $u\in\naturals_0$; solving for $\numeratorBase{v,u+1}\cdot \brackets{\crossSection{u+1}\leftShift v}/\denominator{v}$, along with the identity 
$ 
\madicShift{v,u}= m^{\mExp{v+u}}\madicShift{v, u+1}  + \madicDigit{v,u},
$ yields the inductive step for the desired equalities in the theorem statement.  

This proof makes use of the following identity:   the term
$$
\digitDifference{v,u}  =  l^{\lExp{v}}\madicDigit{v,u}-m^{\mExp{v+u}}\ladicDigit{v+u+1,u}.
$$ When $u=0$, this equality follows from the admissibility criterion
\[
\iterateRemainder{v} = \frac{m^{\mExp{v}}\addend{v}+ \translationValue{v} }{l^{\lExp{v}}}.
\] When $u>0$, this equality follows from  the identity 
\[
 \madicDigit{v,u} = \frac{m^{\mExp{v+u}}\ladicDigit{v+u+1,u} +\madicDigit{v+1,u-1} - \ladicDigit{v+u,u-1} }{l^{\lExp{v}}}
 \]
from Corollary \ref{corollary::madicDiagonalDifferences}.


 The proof is by induction on $u$ for $u\in \naturals_0$. When $u=0$,    the identity  $\denominator{v}\madicShift{v,0} =\numeratorBase{v,0}\cdot \brackets{\crossSection{0}\leftShift v} = \numerator{v}$ follows from (\ref{equation::IterateRationalExpression}) in the previous section.  
Thus, we  proceed with the following chain of equalities:
\begin{alignat}{2}
\denominator{v}\madicShift{v,u}&=\numeratorBase{v,u}\cdot \brackets{\crossSection{u}\leftShift v}   \notag\\
&=\numeratorBase{v,u}\cdot\left( l^{\lExp{v+w}}\madicDigit{v+w,u}-m^{\mExp{v+u+w}}\ladicDigit{v+u+w+1,u}\right)_{w=0}^{\tau-1} \notag\\
&=\numeratorBase{v,u}\cdot\left( l^{\lExp{v+w}}\madicDigit{v+w,u}\right)_{w=0}^{\tau-1} - \numeratorBase{v,u}\cdot\left(m^{\mExp{v+u+w}}\ladicDigit{v+u+w+1,u}\right)_{w=0}^{\tau-1}. \notag
\end{alignat} We assert the equality
\begin{equation}\label{equation::crossSectionInductiveStep}
\numeratorBase{v,u}\cdot\left( l^{\lExp{v+w}}\madicDigit{v+w,u}\right)_{w=0}^{\tau-1} - \denominator{v}\madicDigit{v,u} = \numeratorBase{v,u}\cdot\left( m^{\mExp{v+u+w}}\madicDigit{v+w+1,u}\right)_{w=0}^{\tau-1};
\end{equation} this claim follows from the identity $\meps{v+u,w+1} = \meps{v+u,w} + \mExp{v+u+w}$ (Identity \ref{identities::ExponentSums}\ref{PrefixExponentSumsPerturb}) and the chain of equalities
\begin{alignat}{2}
 \numeratorBase{v,u}\cdot\left( l^{\lExp{v+w}}\madicDigit{v+w,u}\right)_{w=0}^{\tau-1} + m^{\meps{v+\tau,\tau}}\madicDigit{v+\tau,u}
 &= \sum_{0\leq w \leq \tau} m^{\meps{v+u,w}}l^{\less{v,\tau-w}} \madicDigit{v+w,u} \notag\\
&= l^{\less{v,\tau}}\madicDigit{v,u} + \sum_{1\leq w \leq \tau} m^{\meps{v+u,w}}l^{\less{v,\tau-w}} \madicDigit{v+w,u} \notag\\
&= l^{\less{v,\tau}}\madicDigit{v,u} + \sum_{0\leq w < \tau} m^{\meps{v+u,w+1}}l^{\less{v,\tau-1-w}} \madicDigit{v+w+1,u} \notag\\
&= l^{\less{v,\tau}}\madicDigit{v,u} +\numeratorBase{v,u}\cdot\left( m^{\mExp{v+u+w}}\madicDigit{v+w+1,u}\right)_{w=0}^{\tau-1}; \notag
\end{alignat} we can write $\madicDigit{v+\tau,u} = \madicDigit{v, u}$ and $m^{\meps{v+\tau,\tau}} = m^{\meps{v,\tau}}$; thus, subtracting $l^{\less{v,\tau}}\madicDigit{v,u}$ throughout yields (\ref{equation::crossSectionInductiveStep}).

We appeal to the Identity \ref{identities::ExponentSums}\ref{PrefixExponentSumsPerturb} once more to write $\meps{v+u,w} + \mExp{v+u+w} = \mExp{v+u} + \meps{v+u+1,w}$, and we can now express the difference
\begin{alignat}{2}
\denominator{v}\madicShift{v,u}-  \denominator{v}\madicDigit{v,u}&= \numeratorBase{v,u}\cdot\left( m^{\mExp{v+u+w}}\madicDigit{v+w+1,u}\right)_{w=0}^{\tau-1} -  \numeratorBase{v,u}\cdot\left(m^{\mExp{v+u+w}}\ladicDigit{v+u+w+1,u}\right)_{w=0}^{\tau-1}\notag\\
&= m^{\mExp{v+u}} \numeratorBase{v,u+1}\cdot \brackets{\crossSection{u+1}\leftShift v} \notag
\end{alignat}
as required.

\end{proof}

Let us establish an upper bound of $3^{\tau}$ for a potential, periodic iterate value over $\naturals$ for the $3x+1$ Problem.  In this context,   the authors in \cite{BelagaMignotte}  have demonstrated that the maximal iterate $\iterate{\max}$ within a periodic orbit admits the upper bound
\[
\iterate{\max}  < \frac{\pfrac{3}{2}^{\tau-1}}{1 - \frac{3^{\tau}}{2^{\less{\tau}}}}\leq  \tau^C\pfrac{3}{2}^{\tau-1} = o\parentheses{3^{\tau-1}}
\] for some effectively computable constant $C$ (by applying the result in \cite{BakerWustholz}). A recent upper bound on $C$  is available in \cite{Rhin}, in which the author establishes the inequality 
\[
\left|-\less{\tau}\log 2 +\tau \log 3\right| \geq \less{\tau}^{-13.3}
\] (in their notation, we set $u_0 = 0$, $u_1 = -\less{\tau}$, and $u_2 = \tau$);  consequently, assuming $2^{\less{\tau}}>3^{\tau}$, we can write \footnote{We can shed the logarithms:  when $|w| < 1$,  the  power series expansion of $\log (1+w) = \sum_{u\geq 1}(-1)^{u-1}\frac{w^u}{u}$ yields $|\log(1+w)| \leq 2|w|$ when $|w| \leq \frac{1}{2}$. See \cite{Evertse} (Corollary 1.6). } 
\[
1 - \frac{3^{\tau}}{2^{\less{\tau}}} \geq  \frac{\less{\tau}^{-13.3}}{2}.
\] According to \cite{Eliahou}, in a periodic orbit over $\naturals$  of length $\less{\tau}$, the ratio $\frac{\less{\tau}}{\tau}$ satisfies the inequality
\[
\frac{\less{\tau}}{\tau} \leq \lg\parentheses{3+\frac{1}{\iterate{\min}}} \leq 2; 
\]   numerical computation yields
\[
\iterate{\max}< \pfrac{3}{2}^{\tau-1} 2\cdot (2\tau)^{13.3} < 3^{\tau} 
\] when $\tau \geq 103$.  
 
 Thus, if $\iterate{\max}> 3^{\tau}$ and $\iterate{\max}\in \naturals$, then $\tau < 103$. However, the author in \cite{Garner} demonstrates that the length of a non-trivial periodic orbit (excluding $1$) over \naturals\ must satisfy the inequality $2\tau \geq \less{\tau} \geq 35,400$.

   Thus, in the case of the $3x+1$ Problem, the iterate 
$\iterate{v}$ is a positive integer if and only if $$\iterate{v} = \brackets{\numerator{v}\denominator{v}^{-1} } \bmod 3^{\tau} =  \brackets{\numerator{v}\denominator{v}^{-1} } \bmod 2^{\less{\tau}}.$$  This condition is met if and only if,  for all $v\in \nZeroOpenSet{\tau}$,
\begin{enumerate}[i.]
 \item the sequence $\left(\madicQuotient{v,u}\right)_{u\geq 0}$ is non-increasing with $\madicQuotient{v,u} \in \naturals_0$, and 
 \item for $y\in \naturals_0$, the equalities  $\madicQuotient{v,\tau+y} = \ladicQuotient{v,\tau+y} = \digitDifference{v,\tau+y} = 0$ hold.
 \end{enumerate}

For the last result, we establish criteria for the integrality of the iterate values by turning our focus to the prefixes of the $\graded{m}$-adic and \graded{l}-adic expansions, which we define here.
\begin{definition}\label{residueDefinition}  For each $u\in \naturals$ and $v\in \nZeroOpenSet{\tau}$, define
\[
\madicResidue{v,u} = \sum_{0\leq w < u} m^{\meps{v,w}}\madicDigit{v,w} = \brackets{\numerator{v}\denominator{v}^{-1}}_{m^{\meps{u}}}
\] and
\[
\ladicResidue{v,u} = \sum_{0\leq w < u} l^{\less{v,w}}\ladicDigit{v,w} = \brackets{\numerator{v}\denominator{v}^{-1}}_{l^{\less{u}}}
\]  to be the {\sl \graded{m}-residue (of index $u$)} and  {\sl \graded{l}-residue (of index $u$)}, respectively, of $\iterate{v} = \numerator{v}/\denominator{v}.$
\end{definition}

We conclude this article by establishing the identities
	\[
\frac{\numerator{v}}{\denominator{v}} =  \madicResidue{v,\tau} + m^{\meps{v,\tau}}\pfrac{\madicResidue{v,\tau} - \ladicResidue{v,\tau}}{\denominator{v}} = \ladicResidue{v,\tau} + l^{\less{v,\tau}}\pfrac{\madicResidue{v,\tau} - \ladicResidue{v,\tau}}{\denominator{v}};
	\] thus, one can address the question of divisibility of $\numerator{v}$ by $\denominator{v}$ by determining whether the difference $\madicResidue{v,\tau} - \ladicResidue{v,\tau} \in \denominator{v}\integers$.  Furthermore, under the assumption that the ratio is integer-valued, we have $\numerator{v}/\denominator{v} \in \naturals$ if and only if $\frac{\madicResidue{v,\tau} - \ladicResidue{v,\tau}}{\denominator{v}} \in \naturals_0$, as $\madicResidue{v,\tau} \in \nZeroOpenSet{m^{\meps{v,\tau}}}$, and $\ladicResidue{v,\tau} \in \nZeroOpenSet{l^{\less{v,\tau}}}$.

To this end, we require the following lemma.

\begin{lemma}\label{lemma::ladicResidueForwardMap}    Assume  \ref{hypotheses::periodicIterates},  \ref{definition::mlSystemsTerms}, \ref{definition::digitDifference},  and \ref{residueDefinition}. 

For $v\in\nZeroOpenSet{\tau}$, the equality
\[
\frac{m^{\mExp{v}}\ladicResidue{v+1,u} + \translationValue{v}}{l^{\lExp{v}}} = \ladicResidue{v,u} + l^{\less{v,u-1}}\digitDifference{v-u,u}
\]
holds for $u\in \naturals$.

\end{lemma}

\begin{proof}  The proof is by induction on $u$. 
For $u=1$, we have the following chain of equalities as a result of the admissibility criterion in (\ref{equation::admissibilityEquation}):
\begin{alignat}{2}
\frac{m^{\mExp{v}}\ladicResidue{v+1,1} + \translationValue{v}}{l^{\lExp{v}}} &= \frac{m^{\mExp{v}}\ladicDigit{v+1,0} + \translationValue{v}}{l^{\lExp{v}}} \notag\\
&= \madicDigit{v,0} \notag\\
&= \ladicDigit{v,0} +\parentheses{\madicDigit{v,0} - \ladicDigit{v,0}}\notag\\
&= \ladicResidue{v,1} + l^{\less{v,0}}\digitDifference{v-1,1}.\notag
\end{alignat}

Thus, assume the claim for $u\in \naturals$. To demonstrate the claim for $u+1$, we begin by writing
\begin{alignat}{2}
\frac{m^{\mExp{v}}\ladicResidue{v+1,u+1} + \translationValue{v}}{l^{\lExp{v}}} &= \frac{m^{\mExp{v}}\ladicResidue{v+1,u}  +  \translationValue{v}+ m^{\mExp{v}}l^{\less{v+1,u}}\ladicDigit{v+1,u}}{l^{\lExp{v}}}  \notag\\
&= \ladicResidue{v,u} +  l^{\less{v,u-1}}\digitDifference{v-u,u} +m^{\mExp{v}}l^{\less{v,u-1}}\ladicDigit{v+1,u}\notag
\end{alignat} as per the inductive assumption and the identity $\less{v+1,u} = \lExp{v} + \less{v,u-1}$ from (\ref{identities::ExponentSums}.\ref{SuffixExponentSumsPerturb}).
According to Theorem \ref{theorem:diagonalDifferences}, we have the identity
\[
\ladicDigit{v+1,u} = \frac{l^{\lExp{v-u}}\madicDigit{v-u,u} - \digitDifference{v-u,u} }{m^{\mExp{v}}};
\] thus, we can write
\[
m^{\mExp{v}}l^{\less{v,u-1}}\ladicDigit{v+1,u} = l^{\less{v,u-1}}\brackets{l^{\lExp{v-u}}\madicDigit{v-u,u}-\digitDifference{v-u,u}}.
\]  We also have $\less{v,u} = \less{v,u-1} + \lExp{v-u}$ from Identity (\ref{identities::ExponentSums}.\ref{SuffixExponentSumsPerturb});  thus, we have
\begin{alignat}{2}
\frac{m^{\mExp{v}}\ladicResidue{v+1,u+1} + \translationValue{v}}{l^{\lExp{v}}} &=  \ladicResidue{v,u} +  l^{\less{v,u-1}}\brackets{l^{\lExp{v-u}}\madicDigit{v-u,u} }\notag\\
&=  \ladicResidue{v,u+1} +  l^{\less{v,u}}\brackets{\madicDigit{v-u,u} - \ladicDigit{v,u}}\notag\\
&=  \ladicResidue{v,u+1} +  l^{\less{v,u}}\digitDifference{v-1-u,u+1}\notag
\end{alignat}
as required.
\end{proof}

This lemma, along with the following (single-radix) graded polyadic division algorithm, yield the final result of this section.

\begin{algorithm}[H]\label{algorithm::pAdicDivision2}
		\KwIn{An integer $n$}
		\KwIn{An sequence of $\tau$ integers $\List{\gradationSymbol}{0}{\tau-1}$}
		\KwIn{An integer $\numeratorSymbol$}
		\KwIn{An integer $M$}
		\KwIn{An integer $\denominatorSymbol$ with $\gcd(\denominatorSymbol,n)=1$}
		\KwOut{The $\factor{n}$-adic canonical expansion of $\frac{\numeratorSymbol}{\denominatorSymbol}$ with (graded) precision $M$}
		$\mathcal{N}_{0} \leftarrow \numeratorSymbol$\;
		\For{$u: 0 \to M-1$}{
			$\gradationSymbol_u \leftarrow \gradationSymbol_{u\bmod \tau}$\;
			$\madicDigit{u} \leftarrow \denominatorSymbol^{-1}\mathcal{N}_{u}\bmod n^ {\gradationSymbol_{ u}}$\; 
			$\mathcal{N}_{u+1} \leftarrow \frac{\mathcal{N}_{u}- \denominatorSymbol\madicDigit{u}}{n^{\gradationSymbol_{ u}}}$\;
		}
		\Return $\parentheses{\madicDigit{u}}_{u=0}^{M-1}$\;
		\caption{The graded polyadic division algorithm with radix $n$, gradation $\List{\gradationSymbol}{0}{\tau-1}$, and precision $M$}
\end{algorithm}

\begin{theorem}\label{theorem::NumeratorAsDifference} Assume   \ref{hypotheses::periodicIterates},  \ref{definition::mlSystemsTerms}, and \ref{residueDefinition}.

For $v\in \nZeroOpenSet{\tau}$, the equality
\[
 \numerator{v} = {l^{\less{v,\tau}}\madicResidue{v,\tau} - m^{\meps{v,\tau}}\ladicResidue{v,\tau}}
\] holds.
\end{theorem}
\begin{proof}
 Let $N =\numerator{v}$, and let $D = \denominator{v}$.
We will prove the claim by applying the graded polyadic division algorithm above to compute the first $\tau$ $\graded{m}$-adic digits of $\frac{N}{D}$.

For $u\in \nZeroSet{\tau}$, define the sum
\[
S_{v,\tau-u} = \sum_{u\leq w < \tau} m^{\meps{v+u,w-u}}l^{\less{v,\tau-1-w}}\translationValue{v+w}.
\] We will  prove the equality 
\[
\madicQuotient{v,u}D = -l^{\less{v,\tau-u}}\ladicResidue{v+u,u} + S_{v,\tau-u} + m^{\meps{v+u,\tau-u}}\madicResidue{v,u}.
\] by induction on $u$;  afterwards, we will derive the desired result by evaluating this equality at $u = \tau$.

For $u=0$, we have 
\[
\madicQuotient{v,0}D = -0+ \sum_{0\leq w < \tau} m^{\meps{v,w}}l^{\less{v,\tau-1-w}}\translationValue{v+w} + 0
\] as required. 

 Assuming the inductive claim for $u\geq 0$, we will write
\begin{alignat}{2}
\madicQuotient{v, u+1}D &= \frac{\madicQuotient{v,u}D - \madicDigit{v,u}D}{m^{\mExp{v+u}}}\notag \\
&= \frac{-l^{\less{v,\tau-u}}\ladicResidue{v+u, u} + S_{v,\tau-u}+ m^{\meps{v+u,\tau-u}}\madicResidue{v,u} - \madicDigit{v,u}\parentheses{l^{\less{v,\tau}} - m^{\meps{v,\tau}}} }{m^{\mExp{v+u}}}.
\label{equation::beforeSubstitution}
\end{alignat}

We will write
\begin{alignat}{2}
S_{v,\tau-u} &= l^{\less{v,\tau-1-u}}\translationValue{v,u} +  \sum_{u+1\leq w < \tau} m^{\mExp{v+u}+\meps{v+u+1,w-u-1}}l^{\less{v,\tau-1-w}}\translationValue{v+w}\notag\\
&= l^{\less{v,\tau-1-u}}\translationValue{v,u} + m^{\mExp{v+u}}S_{v,\tau-u-1}.\notag
\end{alignat}

According Identity \ref{identities::ExponentSums} on exponent sums, we have the equality
$
 \meps{v,\tau} =  \meps{v,u} + \meps{v+u,\tau-u};
$
thus, we will also write
\begin{alignat}{2}
m^{\meps{v+u,\tau-u}}\madicResidue{v, u} +m^{\meps{v,\tau}}\madicDigit{v,u} &= m^{\meps{v+u,\tau-u}}\parentheses{\madicResidue{v, u} + m^{\meps{v,u}}\madicDigit{v,u}}\notag\\
&= m^{\mExp{v+u} + \meps{v+u+1,\tau-u-1}  }\madicResidue{v, u+1}.\notag
\end{alignat}

According to Lemma \ref{lemma::ladicResidueForwardMap}, the expression
\[
-l^{\less{v,\tau-u}}\ladicResidue{v+u, u} - l^{\less{v,\tau}}\madicDigit{v,u}+ l^{\less{v,\tau-1-u}}\translationValue{v,u} =- m^{\mExp{v+u}}l^{\less{v,\tau-u-1}}\ladicResidue{v+u+1, u+1}.
\] 
Thus, we can substitute these expressions in (\ref{equation::beforeSubstitution}) above and write
\[
\madicQuotient{v, u+1}D =  -l^{\less{v,\tau-u-1}}\ladicResidue{v+u+1 , u+1} + S_{v,\tau-u-1} + m^{\meps{v+u+1,\tau-u-1}}\madicResidue{v,u+1}
\]
as required.


We can write $\ladicResidue{v+\tau,\tau} = \ladicResidue{v,\tau}$; thus, when $u = \tau$,  we have the equality
\begin{equation}\label{equation::tauIterateShift}
\madicQuotient{v,\tau}D  = \madicResidue{v,\tau} - \ladicResidue{v,\tau}.
\end{equation}We now can write
\begin{equation}\label{equation::madicRepresentationOfIterate}
\frac{N}{D} =  \madicResidue{v,\tau} + m^{\meps{v,\tau}}\madicShift{v,\tau}= \madicResidue{v,\tau} + m^{\meps{v,\tau}}\pfrac{\madicResidue{v,\tau} - \ladicResidue{v,\tau}}{D};
\end{equation}  when  solving for $N$, we arrive at the equalities
\[
N  = D\madicResidue{v,\tau} + m^{\meps{v,\tau}}\parentheses{\madicResidue{v,\tau} - \ladicResidue{v,\tau}} = {l^{\less{v,\tau}}\madicResidue{v,\tau} - m^{\meps{v,\tau}}\ladicResidue{v,\tau}}
\] as required.

\end{proof}

We have the following consequence of Theorem \ref{theorem::NumeratorAsDifference}. 

\begin{corollary}\label{corollary::iterateResidueRepresentation}
For each $v\in \nZeroOpenSet{\tau}$, the iterate $\frac{\numerator{v}}{\denominator{v}}$ satisfies the chain of equalities $$\frac{\numerator{v}}{\denominator{v}} = \madicResidue{v,\tau} + m^{\meps{v,\tau}}\pfrac{\madicResidue{v,\tau}  - \ladicResidue{v,\tau} }
{\denominator{v}} = \ladicResidue{v,\tau} + l^{\less{v,\tau}}\pfrac{\madicResidue{v,\tau}  - \ladicResidue{v,\tau} }
{\denominator{v}}.$$
\end{corollary}

\begin{proof}
The first equality is a result of (\ref{equation::tauIterateShift}) yielding (\ref{equation::madicRepresentationOfIterate}).  As per Corollary \ref{corollaryQuotientIdentity}, we have the chain of equalities $\madicQuotient{v,\tau}=\ladicQuotient{v+\tau,\tau} = \ladicQuotient{v,\tau}$; thus, we  can recast (\ref{equation::tauIterateShift}) as $\ladicQuotient{v,\tau} = \frac{\madicResidue{v,\tau} - \ladicResidue{v,\tau}}{\denominator{v}}$ to derive the second equality.
\end{proof}

As $m$ and $l$ are coprime, we can write
\[
1 = l^{\less{\tau}}\bracketsInverse{l^{\less{v,\tau}}}_{m^{\meps{v,\tau}}} - m^{\meps{v,\tau}}\bracketsInverse{-m^{\meps{v,\tau}}}_{l^{\less{v,\tau}}} =  l^{\less{v,\tau}}\bracketsInverse{\denominator{v}}_{m^{\meps{v,\tau}}} - m^{\meps{v,\tau}}\bracketsInverse{\denominator{v}}_{l^{\less{v,\tau}}}
\] as per B\'{e}zout's Identity. As per Theorem \ref{theorem::NumeratorAsDifference}, we can write
\[
\numerator{v}  =  l^{\less{v,\tau}}\bracketsInverse{\numerator{v}\denominator{v}}_{m^{\meps{v,\tau}}} - m^{\meps{v,\tau}}\bracketsInverse{\numerator{v}\denominator{v}}_{l^{\less{v,\tau}}}.
\]

Another consequence of Theorem \ref{theorem::NumeratorAsDifference}:
Fix $\tau,\meps{\tau},\less{\tau} \in \naturals$ where both $\meps{\tau}\geq \tau$ and $\less{\tau} \geq \tau$.  Let $\graded{f} = (\mExp{w})_{w=0}^{\tau-1}$ and  $\graded{e} = (\lExp{w})_{w=0}^{\tau-1}$ be  positive partitions of order $\tau$ of $\meps{\tau}$ and $\less{\tau}$, respectively.
For any pair $(m,l)$ of coprime natural numbers greater than or equal to $2$,  we can construct a $(m,l)$-system  that admits a sequence with the translation values $\graded{a} = \graded{1}^{\tau}$: we prescribe $\iterateRemainder{w} \equivMod{m^{\mExp{w}}}\brackets{l^{\lExp{w}}}^{-1}$ and $\addend{w} \equivMod{l^{\lExp{w}}} \brackets{-m^{\mExp{w}}}^{-1}$ for each $w\in \nZeroOpenSet{\tau}$.

By expressing a periodic iterate within such a system as $\iterate{0} = \frac{\numeratorBase{0}\cdot{\graded{a}}}{\denominator{}}$,  we can write
\[
\numeratorBase{0}\cdot{\graded{a}} =  \sum_{0\leq w < \tau}m^{\meps{w}}l^{\less{\tau-1-w}} =l^{\less{\tau}}\madicResidue{\tau} -  m^{\meps{\tau}}\ladicResidue{\tau}.
\] When the terms $m=3$ and $l=2$, we have a {\sl three-smooth representation} (\cite{Avidon}, \cite{Erdos}, \cite{Selfridge}) of the numerator of the iterate.

\bibliographystyle{plain}

\end{document}